\documentclass{amsart}%
\usepackage{amsmath}
\usepackage{amssymb}
\usepackage{graphicx}
\usepackage{amscd}
\usepackage{amsfonts}%
\providecommand{\U}[1]{\protect\rule{.1in}{.1in}}
\newtheorem{theorem}{Theorem}
\theoremstyle{plain}

\newtheorem{corollary}{Corollary}

\newtheorem{definition}{Definition}

\newtheorem{lemma}{Lemma}

\newtheorem{remark}{Remark}

\numberwithin{equation}{section}

\begin{document}
\title[\textsf{The Structure of NAFIL Loops}]{\textsf{The Structure of Non-Associative Finite Invertible Loops (NAFIL)*}}
\author{\textbf{Raoul E. Cawagas}}
\address{Raoul E. Cawagas, SciTech R\&D Center, Polytechnic University of the
Philippines, Sta. Mesa, Manila}
\email{raoulec1@yahoo.com}
\thanks{*This paper is part of a series of papers on Non-Associative Finite Invertible
Loops (NAFIL).}

\begin{abstract}
The NAFIL is a loop in which every element has a unique (two-sided) inverse.
NAFIL loops can be classified into two types: \emph{composite} (with at least
one non-trivial subsystem) and \emph{non-composite} or \emph{plain} (without
any non-trivial subsystem). This paper deals with the structure of these two
types of loops. In particular we shall introduce an important class of
composite loops called \emph{block products}.

\end{abstract}
\maketitle

\section{\textbf{Introduction}}

In a previous paper [2] we introduced the concept of a finite algebra and
defined the \emph{non-associative finite invertible loop (NAFIL)}. This is a
loop in which every element has a unique (two-sided) inverse and it satisfies
all group axioms except the Associative axiom. Because of this, the NAFIL has
many structural features in common with the group. However, because the NAFIL
is non-associative, it has its own distinctive structure that sets it apart
from most algebraic systems.

The study of NAFIL structure is in its early stages and not much information
about it is available in the current literature (English). In this paper, we
shall begin the study of NAFIL structure by considering the basic properties
of NAFIL loops with non-trivial subsystems (\textbf{\emph{composite}}), and
those without any non-trivial subsystems (\textbf{\emph{non-composite }%
}or\textbf{\emph{\ plain}}).

Before we proceed, let us recall the idea of an \emph{abstract mathematical
system} [2, 5]. Such a system essentially consists of: a non-empty set $S$ of
distinct elements, at least one binary operation $\star$, an equivalence
relation $=$, a set of axioms, as well as a set of definitions and theorems
and is usually denoted by $\left(  S,\star\right)  $. The heart of the system
is the set of axioms from which all the theorems are derived.

Most algebraic systems like groups, rings, and fields satisfy some or all of
the following axioms or postulates:

\begin{itemize}
\item A1 - For all $a,b\in S,$ $a\star b\in S$. \emph{(Closure axiom)}

\item A2 -There exists a unique element $e\in S,$ called the \emph{identity},
such that $e\star a=a\star e=a$ for all $a\in S$. \emph{(Identity axiom)}

\item A3 - Given an identity element $e\in S,$ for every $a\in S$ there exists
a unique element $a^{-1}\in S,$ called its \emph{inverse}, such that $a\star
a^{-1}=a^{-1}\star a=e$. \emph{(Inverse axiom)}

\item A4 - For every $a,b\in$ $S$ there exists unique $x$, $y\in$ $S$ such
that $a\star x=b$ and $y\star a=b$. \emph{(Unique Solution axiom)}

\item A5 - For every $a,b\in S$, $a\star b=b\star a$\emph{. (Commutative
axiom)}

\item A6 - For all $a,b,c\in S$, $(a\star b)\star c=a\star(b\star c)$.
\emph{(Associative axiom)}\medskip
\end{itemize}

If the set $S$ is of finite order, then $(S,\star)$ is called a \emph{finite
system}. In this paper, we shall only consider finite systems.\medskip

The simplest algebraic system is the \emph{groupoid}; it is only required to
satisfy A1. This is a trivial system and not much can be said about it. A
groupoid that also satisfies A4 is called a \emph{quasigroup}; a quasigroup
that satisfies A2 is called a \emph{loop}; and a loop that satisfies A3 is
called an \emph{invertible loop}. Moreover, an invertible loop that satisfies
A6 is called a \emph{group} and one that does not satisfy A6 is called a
\emph{NAFIL} (\emph{non-associative finite invertible loop}). In this paper,
the word \emph{groupoid} will often be used as a generic term for any system
satisfying A1 (like quasigroups, loops, and groups).

\section{\textbf{Composite NAFIL Loops}}

Most mathematical systems contain smaller systems in their structures. Thus,
sets can have subsets and groups can have subgroups. Similarly, NAFIL loops
can also have subsystems (like loops or groups). To be specific, a subsystem
that is a group shall be called a \textbf{\emph{subgroup}} while one that is a
loop (like the NAFIL) shall be called a \textbf{\emph{subloop}} [6].
Otherwise, the word \textbf{\emph{subsystem}} shall be used as a generic term
for both \emph{subgroup} and \emph{subloop}.

\begin{definition}
Let $\left(  S,\star\right)  $ be a loop of order $n$, and let $\overline{S}$
be a non-empty subset of order $m$ of $S$. If $\left(  \overline{S}%
,\star\right)  $ satisfies all axioms satisfied by $\left(  S,\star\right)  $,
then it is called a \textbf{subsystem} of $\left(  S,\star\right)  $.
Moreover, a subsystem of order $m$ is also called \textbf{proper} if $1\leq
m<n$ and \textbf{improper} if $m=n$. If $m=1,$ the proper subsystem is called
\textbf{trivial}; otherwise it is called \textbf{non-trivial}.
\end{definition}

This definition holds for any groupoid in general. It simply states that if a
given system satisfies a set of axioms, then any of its subsystems must also
satisfy this set of axioms completely. Thus, a subgroup (or subloop) must
satisfy all group (or loop) axioms. Since the group satisfies all NAFIL
axioms, then a NAFIL can have a group as a subsystem.

A loop $\left(  \mathcal{L},\star\right)  $ of order $n$ will be called
\textbf{\emph{composite}} if it contains at least one non-trivial proper
subsystem. A loop is called \textbf{\emph{Lagrangian}} if the order $m$ of any
of its subsystems is always a divisor (or factor) of $n$. Otherwise, it is
called \textbf{\emph{non-Lagrangian}. }Moreover\textbf{, }a subloop $\left(
H,\star\right)  $ of a loop $\left(  \mathcal{L},\star\right)  $ is also said
to be \textbf{\emph{Lagrange-like}} [6] if the order of $H$ divides the order
of $\mathcal{L}.$

Following common practice, we shall also denote the order $n$ of any set $S$
by $\left|  S\right|  .$ Moreover, we shall write $\overline{S}\leq S$ to mean
that $\overline{S}$ is a subsystem of $S.$

In what follows, we shall introduce and prove a lemma and a number of
interesting theorems on subsystems of NAFIL loops. For convenience, we shall
often use $\mathcal{L}$ to denote the set of elements $\{\ell_{i}\mid
i=1,...,n\}$ of any finite invertible loop, $\overline{\mathcal{L}}$ to denote
a subset of $\mathcal{L}$, $1\equiv\ell_{1}$ to denote the identity element,
and $\ell_{i}\ell_{j}$ to denote the product $\ell_{i}\star\ell_{j},$ that is,
$\ell_{i}\ell_{j}\equiv$ $\ell_{i}\star\ell_{j}$ if no confusion arises.

\begin{lemma}
Let $\left(  \mathcal{L},\star\right)  $ be an invertible loop of order $n$
and let $\overline{\mathcal{L}}$ be a non-empty subset of order $m$ of
$\mathcal{L}$. If $\overline{\mathcal{L}}$ is closed under $\star$ and
$T=\{\ell_{y}\in\mathcal{L}\mid\ell_{y}\notin\overline{\mathcal{L}}\}$, then:
$\ell_{x}\ell_{y}$ and $\ell_{y}\ell_{x}$ are in $T$ for all $\ell_{x}%
\in\overline{\mathcal{L}}$ and all $\ell_{y}\in T$.
\end{lemma}

%

\proof
Let $\overline{\mathcal{L}}=\{\ell_{1},...,\ell_{m}\}$ and $T=\{\ell
_{m+1},...,\ell_{n}\}.$ Since $T=\{\ell_{y}\in\mathcal{L}\mid\ell_{y}%
\notin\overline{\mathcal{L}}\}$, then it follows that $\mathcal{L}%
=\overline{\mathcal{L}}\cup T$ and $\overline{\mathcal{L}}\cap T=\emptyset$.
Hence,
\begin{equation}
\ell_{1},...,\ell_{m},\ell_{m+1},...,\ell_{n} \tag{L1.1}%
\end{equation}
are the $n$ distinct elements of $\mathcal{L}$. Now, let $\ell_{x}$ be any
element of $\mathcal{L}$ and form the sets:
\begin{equation}
\ell_{x}\overline{\mathcal{L}}=\{\ell_{x}\ell_{1},...,\ell_{x}\ell
_{m}\}\text{\quad and\quad}\ell_{x}T=\{\ell_{x}\ell_{m+1},...,\ell_{x}\ell
_{n}\} \tag{L1.2}%
\end{equation}%
\begin{equation}
\overline{\mathcal{L}}\ell_{x}=\{\ell_{1}\ell_{x},...,\ell_{m}\ell_{x}%
\}\quad\text{and\quad}T\ell_{x}=\{\ell_{m+1}\ell_{x},...,\ell_{n}\ell_{x}\}
\tag{L1.3}%
\end{equation}
Since $\left(  \mathcal{L},\star\right)  $ is an invertible loop and
$\overline{\mathcal{L}}\cap T=\emptyset$, it follows from A1, A4 and Eqs.
(L1.2) and (L1.3) that $\ell_{x}\mathcal{L}=(\ell_{x}\overline{\mathcal{L}%
})\cup(\ell_{x}T)=\mathcal{L}$ and $\mathcal{L}\ell_{x}=(\overline
{\mathcal{L}}\ell_{x})\cup(T\ell_{x})=\mathcal{L}.$ Next, let $\ell_{x}%
\in\overline{\mathcal{L}}$. Since $\overline{\mathcal{L}}$ is closed under
$\star$, then $\ell_{x}\overline{\mathcal{L}}=\overline{\mathcal{L}}\ell
_{x}=\overline{\mathcal{L}}$ and hence it follows that $\ell_{x}T=T\ell
_{x}=T.$ This implies that if $\ell_{x}\in\overline{\mathcal{L}}$ and
$\ell_{y}\in T$, then $\ell_{x}\ell_{y}$, $\ell_{y}\ell_{x}\in T$. This
completes the proof of the Lemma. $\blacksquare$\vspace{0.05in}

With this Lemma, we can now prove a number of interesting theorems on
subsystems of NAFIL loops.

\begin{theorem}
Let $\left(  \mathcal{L},\star\right)  $ be a NAFIL of order $n$ and let
$\overline{\mathcal{L}}$ be a non-empty subset of order $m$ of $\mathcal{L}$.
If $\star$ is closed on $\overline{\mathcal{L}}$, then $\left(  \overline
{\mathcal{L}},\star\right)  $ is a subsystem of $\left(  \mathcal{L}%
,\star\right)  $.
\end{theorem}

%

\proof
By hypothesis, $\star$ is closed on $\overline{\mathcal{L}}$. Hence, $\left(
\overline{\mathcal{L}},\star\right)  $ satisfies A1. By Definition 1 if
$\left(  \overline{\mathcal{L}},\star\right)  $ is a subsystem of $\left(
\mathcal{L},\star\right)  $, then it must also satisfy A2, A3, and A4. We
first prove that $\left(  \overline{\mathcal{L}},\star\right)  $ satisfies A4.
Let $\ell_{a},\ell_{b}\in\overline{\mathcal{L}}$ and let $T=\{\ell_{r}%
\in\mathcal{L}\mid\ell_{r}\notin\overline{\mathcal{L}}\}$. Since $\ell
_{a},\ell_{b}$ are also elements of $\mathcal{L}$, there exist unique elements
$\ell_{x}$ and $\ell_{y}$ such that $\ell_{a}\ell_{x}=\ell_{b}$ and $\ell
_{y}\ell_{a}=\ell_{b}$. Hence, either $\ell_{x}\in\overline{\mathcal{L}}$ or
else $\ell_{x}\in T$. Suppose $\ell_{x}\in T$. Then by Lemma 1, $\ell_{a}%
\ell_{x}=\ell_{b}\in T$. This is a contradiction since $\ell_{b}\in
\overline{\mathcal{L}}$. Therefore, $\ell_{x}\in\overline{\mathcal{L}}$.
Likewise, it follows from Lemma 1 that $\ell_{y}\in\overline{\mathcal{L}}$.
Hence, A4 is satisfied and $\left(  \overline{\mathcal{L}},\star\right)  $ is
at least a quasigroup.

Let $\ell_{1}$ be the identity element of $\mathcal{L}$. Assume that $\ell
_{1}\in T$ and let $\ell_{x}\in\overline{\mathcal{L}}$. Then, by Lemma 1,
$\ell_{x}\ell_{1}=\ell_{x}\in T$ which is false since, by hypothesis,
$\ell_{x}\in\overline{\mathcal{L}}$. This implies that $\ell_{1}\in
\overline{\mathcal{L}}$ and hence A2 is satisfied.

To show that every element in $\overline{\mathcal{L}}$ has an inverse in
$\overline{\mathcal{L}}$, let $\ell_{a}\in\overline{\mathcal{L}}$. Since
$\ell_{1}\in\overline{\mathcal{L}}$, there exists $\ell_{x}$, $\ell_{y}%
\in\overline{\mathcal{L}}$ satisfying: $\ell_{a}\ell_{x}=\ell_{1}$ and
$\ell_{y}\ell_{a}=\ell_{1}$. By A4, the solution to the these equations is
$\ell_{x}=\ell_{y}=\ell_{a}^{-1}$, where $\ell_{a}^{-1}$ is the inverse of
$\ell_{a}$ in $\mathcal{L}$ under $\star$. Hence, $\ell_{a}^{-1}$ is in
$\overline{\mathcal{L}}$ and A3 is satisfied. Therefore, $\left(
\overline{\mathcal{L}},\star\right)  $ satisfies A1,A2,A3, and A4 and is thus
a subsystem of $\left(  \mathcal{L},\star\right)  $. This completes the proof
of Theorem 1. $\blacksquare$\vspace{0.05in}

\begin{theorem}
The order $m$ of any proper subsystem of a composite NAFIL of order $n$ is
equal to or less than $n/2$, that is, $m\leq n/2.$
\end{theorem}

%

\proof
Let $\left(  \mathcal{L},\star\right)  $ be a composite NAFIL of order $n$ and
let $\left(  \overline{\mathcal{L}},\star\right)  $ be any proper subsystem of
order $m$. Let $\overline{\mathcal{L}}=\{\ell_{1},....,\ell_{m}\}$ and
$T=\{\ell_{m+1},...,\ell_{n}\}$ such that $T=\{\ell_{j}\in\mathcal{L}\mid
\ell_{j}\notin\overline{\mathcal{L}}\}$. Then $\ell_{1},...,\ell_{m}%
,\ell_{m+1},...,\ell_{n}$ are the $n$ distinct elements of $\mathcal{L}$. Now,
let $\ell_{j}$ be any element of $T$. Form the $n$ products:
\[
\ell_{j}\ell_{1},...,\ell_{j}\ell_{m},\,\ell_{j}\ell_{m+1},...,\ell_{j}%
\ell_{n}%
\]
By A1 and A4, these $n$ products are the $n$ distinct elements of
$\mathcal{L}$ in some order. By Lemma 1, if $\ell_{i}\in\overline{\mathcal{L}%
}$ and $\ell_{j}\in T$, then $\ell_{j}\ell_{i}\in T$. Therefore, the $m$
distinct products $\ell_{j}\ell_{1},...,\ell_{j}\ell_{m}$ are all elements of
$T$. This means that $T$ contains at least these $m$ distinct elements. Since
$T $ obviously contains $n-m$ distinct elements, then the following relation
must hold: $(n-m)\geq m$. This numerical inequality implies that $m\leq n/2.$
Since $m$ also corresponds to the order of $\overline{\mathcal{L}}$, this
proves Theorem 2. $\blacksquare$\vspace{0.05in}

This theorem indicates that unlike groups certain NAFIL loops of order $n$ can
have a subsystem whose order $m$ is not a divisor of $n$. Such a system is,
therefore, non-Lagrangian\textbf{.} Indeed, there are many known NAFIL loops
of this type. Thus, the smallest NAFIL loop $(L_{5},\ast)$ of order $n=5$ is
non-Lagrangian because it has subsystems of order $n=2$.

\begin{center}
$%
\begin{tabular}
[c]{|c|ccccc|}\hline
$\star$ & 1 & 2 & 3 & 4 & 5\\\hline
1 & 1 & 2 & 3 & 4 & 5\\
2 & 2 & 1 & 5 & 3 & 4\\
3 & 3 & 4 & 1 & 5 & 2\\
4 & 4 & 5 & 2 & 1 & 3\\
5 & 5 & 3 & 4 & 2 & 1\\\hline
\end{tabular}
\medskip$

$\underset{\text{Table 1. Cayley table of NAFIL loop }(L_{5},\ast)\text{ of
order n = 5\bigskip}}{}$
\end{center}

\begin{remark}
\emph{Lemma 1 and Theorems 1 and 2 hold for loops in general. In fact, they
hold even for quasigroups since their proofs essentially depend on A1 and
A4.}\medskip
\end{remark}

\subsection{\textbf{Lagrangian Systems}}

Composite systems play a central role in the study of NAFIL structure. One of
the most important classes of\ composite NAFIL loops is the class of
\textbf{\emph{block product systems}} [3] that includes \emph{direct products}
and \emph{coset products} which are Lagrangian.

The idea of the \emph{block product} is a generalization of the \emph{direct
product} concept and is related to groups with a \emph{factor group }(or
\emph{group of cosets}) [9] in group theory. This arises from the observation
that the Cayley table of a group with a normal subgroup, when its entries are
arranged in terms of the cosets of this subgroup, is seen to split up into
blocks that is induced by the group operation on the cosets. The entries in
each block (called a \textbf{\emph{coset block}}) all belong to a single coset
so that each coset can be considered as a single element. These cosets give
rise to a partition of the group elements and they form, under certain
conditions, a group called the factor group.

The concept of the \textbf{\emph{coset }}[8, 9] in finite group theory is of
great importance in the study of associative algebraic structures. The proof
of Lagrange's theorem (that the order of a subgroup is a divisor of the order
of the group) is based on cosets. This concept, however, does not depend on
the associative axiom A6 and it also applies to NAFIL loops and loops in general.

\begin{definition}
Let $\left(  H,\star\right)  $ be a subsystem of a composite system $\left(
\mathcal{L},\star\right)  $ and let $a\in\mathcal{L}$. The subsets
$aH=\{a\star\ell\mid\ell\in H\}$ and $Ha=\{\ell\star a\mid\ell\in H\}$ are
called the \textbf{left} and \textbf{right cosets} of $H$ in $\left(
\mathcal{L},\star\right)  ,$ respectively. Here, the element $a$ is called a
\textbf{coset} \textbf{representative}.
\end{definition}

If $K\subseteq\mathcal{L}$ and $K=aH$ $\mathbf{(}$or $K=Ha\mathbf{)}$, then we
also say that $K$ is a left (or right) coset \textbf{\emph{modulo}} $H$ for
some $a\in\mathcal{L}$.

Not all composite groups, however, have factor groups (or coset groups). In
studying them, the problem posed is to determine the precise conditions under
which the elements of such a group can split up into coset blocks with a
well-defined operation induced by the group operation on the cosets. A
sufficient condition for this is that every left coset is also a right coset.
This condition, however, does not hold in general for loops.

\subsubsection{\textbf{The Block Product}}

In this section, we shall present an important system analogous to groups with
factor groups that applies to both loops and groups.

\begin{definition}
Let $S=\left\{  s_{1},\ldots,s_{n}\right\}  $ be a set of order $n=km$ and let
$B=\left\{  B_{1},\ldots,B_{k}\right\}  $ be a partition of $S$ where every
$B_{i}\in B$ is of order $m.$ Let $\times$ and $\diamond$ be quasigroup-type
operations on the sets $B$ and $S,$ respectively, such that: If $s_{i}\in
B_{p},$ $s_{j}\in B_{q}$ and $B_{p}\times B_{q}=B_{r},$ then $s_{i}\diamond
s_{j}=s_{h}$ for some $s_{h}\in B_{r}.$ The operations $\times$ and $\diamond$
give rise to two quasigroup-type systems $(B,\times)$ called a \textbf{factor
system} and $(S,\diamond)$ called a \textbf{block product}.
\end{definition}

In this definition the block product $\left(  S,\diamond\right)  $, by
analogy, corresponds to the group with a factor group and the factor system
$(B,\times)$ corresponds to the factor group (where the cells of the partition
$B$ take the place of the cosets and the operation $\times$ is called
\textsl{cell multiplication}.). However, the roles of the operations
$\diamond$ and $\times$ can be viewed in two ways: \emph{(a)} the operation
$\times$ is induced by the operation $\diamond$ on $B$ as in the case of
factor groups and \emph{(b)} $\diamond$ is induced by $\times$ on $S.$ In
either case, the definition does not completely specify these quasigroup
operations; it only states a necessary condition that $\diamond$ and $\times$
must satisfy. This incompleteness is what makes the block product concept a
useful tool in the construction of composite algebraic systems. It allows us
to impose certain requirements on the operations $\diamond$ and $\times$ as
well as on the sets $S$ and $B$ to obtain the desired block product.%

\[
\mathbf{The\ Multi-}\Phi\mathbf{\ System}%
\]
\textbf{\ }

Definition 3 does not explicitly define the block product which can assume
several forms. In what follows we shall therefore introduce a particular type
of the block product in line with case (a) in which $\times$ is induced by
$\diamond$ on $B.$ To do this, we need to present first the concept of the
\emph{multi-}$\phi$ \emph{system}.

\begin{definition}
\textbf{\ }Let $C=\{c_{1},...,c_{m}\}$ be a set of $m$ elements and let
$\Phi=\{\phi_{1},...,\phi_{g}\}$ be a set of $g$ closed binary operations on
the set $C.$ The system $\left(  C,\Phi\right)  $ is called a\textbf{\ multi-}%
$\phi$\textbf{\ system} of order $(m;g)$ if it satisfies the following: (I)
The system $\left(  C,\Phi\right)  $ is at least a quasigroup under every
operation $\phi_{x}\in\Phi,$ and (II) Two binary operations $\phi_{u}$%
,$\phi_{v}\in\Phi$ are \textbf{equal}, that is $\phi_{u}=\phi_{v}$, if and
only if $c_{i}\phi_{u}c_{j}=c_{i}\phi_{v}c_{j}$ for all $c_{i},c_{j}\in C$. If
$\left(  C,\Phi\right)  $ is of order $(m;1)$, then it is an ordinary finite
system called a \textbf{mono-}$\phi.$
\end{definition}

The multi-$\phi$ system $\left(  C,\Phi\right)  $ consists of a number of
systems $\left(  C,\phi_{x}\right)  $ with a common set of elements $C.$ This
system is equivalent to what is known as an \emph{indexed algebra}%
\textbf{\emph{\ }} By indexing the operations $\phi_{x}\in\Phi,$ it is
possible to compare or distinguish the operations $\phi_{x}$ and $\phi_{y}$ of
any two systems under $\left(  C,\Phi\right)  $ by means of their indices $x$
and $y.$

Because the multi-$\phi$ system involves operations of various kinds, it
becomes necessary to classify systems into types according to the axioms that
they are required to satisfy (or not satisfy).

Depending on what axioms a system $\left(  S,\star\right)  $ is required to
satisfy, it is usually classified in terms of its \textsl{axiom type} as
follows:\bigskip

\begin{center}
$\underset{\text{Table 2. Axiom Types of some algebraic systems. (Note:
}\emph{(}\sim\emph{A6)}\text{ means \emph{non-associative}.)\bigskip}%
}{\underset{}{
\begin{tabular}
[c]{|l|l|l|}\hline
SYSTEM NAME & AXIOMS SATISFIED & AXIOM TYPE\\\hline%
\begin{tabular}
[c]{l}%
Groupoid\\
Quasigroup\\
Loop\\
Invertible loop\\
NAFIL loop\\
Semigroup\\
Monoid\\
Group\\
Abelian
\end{tabular}
&
\begin{tabular}
[c]{l}%
\emph{A1}\\
\emph{A1, A4}\\
\emph{A1, A4, A2}\\
\emph{A1, A4, A2, A3}\\
\emph{A1, A4, A2, A3; (}$\sim$\emph{A6)}\\
\emph{A1, A6}\\
\emph{A1, A6, A2}\\
\emph{A1, A4, A2, A3, A6}\\
Any system satisfying \emph{A5}%
\end{tabular}
&
\begin{tabular}
[c]{l}%
\emph{A[1]}\\
\emph{A[1,4]}\\
\emph{A[1,4,2]}\\
\emph{A[1,4,2,3]}\\
\emph{A[1,4,2,3](}$\sim$\emph{A6)}\\
\emph{A[1,6]}\\
\emph{A[1,6,2]}\\
\emph{A[1,4,2,3,6]}\\
\emph{A[5]}%
\end{tabular}
\\\hline
\end{tabular}
}}$
\end{center}

A system $\left(  S,\ast\right)  $ is of \textsl{axiom type} \emph{A[1,...,}%
x\emph{]} if it satisfies axioms \emph{A1,...,A}x. This does not imply that
these are the only axioms that $\left(  S,\ast\right)  $ satisfies; it only
means that they are specifically required to be satisfied. Thus, an invertible
loop is of type \emph{A[1,4,2,3] }although it may happen that it also
satisfies \emph{A5} and is therefore abelian. If a system is required not to
satisfy a given axiom \emph{A}x$,$ we shall indicate this by writing
\emph{(}$\thicksim$\emph{A}x\emph{)}$.$ Thus, the axiom type of the NAFIL is
written as \emph{A[1,4,2,3](}$\sim$\emph{A6).}

Henceforth, we shall consider only systems of at least type \emph{A[1,4] }to
be called \textbf{\emph{quasigroup-type}} systems. This includes
\emph{quasigroups}, \emph{loops}, \emph{NAFIL loops}, and \emph{groups}. Note
that the group and the NAFIL satisfy all \emph{invertible loop} axioms.
Because of this, the \emph{group} can also be considered as an
\emph{associative invertible loop} while the NAFIL is a \emph{non-associative
invertible loop. }Thus, the term\emph{\ invertible loop \ }is a generic term
for both groups and NAFIL loops.

Since $\left(  C,\Phi\right)  $ is a multi-$\phi$ system, it represents
several systems, that is, every $\phi_{f}\in\Phi$ defines a system $\left(
C,\phi_{f}\right)  $ of a given axiom type \emph{A[1,...,}x\emph{]}. If these
systems are all of the same axiom type \emph{A[1,...,}x\emph{],} then we
simply call $\left(  C,\Phi\right)  $ by this common axiom type. Thus, if
every $\phi_{f}\in\Phi$ is a quasigroup operation, then we shall call $\left(
C,\Phi\right)  $ a \emph{quasigroup system} for convenience.\bigskip

\subsubsection{\textbf{The Block Product as a Generalized Direct Product}}

We can now introduce an equivalent definition of the block product given by
Definition 3. This is a generalization of the \emph{direct product} [5, 10] in
group theory that involves the \emph{multi}-$\phi$ system.

\begin{definition}
Let $\left(  E,\ast\right)  $ and $\left(  C,\Phi\right)  $ be two
quasigroup-type systems of orders $k$ and $m,$ respectively, where
$\Phi=\left\{  \phi_{pq}\mid p,q=1,...,k\right\}  $ is a set of $k^{2}$
quasigroup operations, and let $\mathcal{L}$ $=E$\textsf{X}$C=\left\{
(e_{i},c_{j})\mid e_{i}\in E,\,c_{j}\in C\right\}  .$ The \textbf{block
product} of $\left(  E,\ast\right)  $ and $\left(  C,\Phi\right)  $ is the
system
\[
\left(  \mathcal{L},\diamond\right)  =\left(  E,\ast\right)  \text{\textsf{X}%
}\left(  C,\Phi\right)
\]
of order $n=km,$ where $\diamond$ is defined by the composition rule:
\begin{equation}
(e_{p},c_{a})\diamond(e_{q},c_{b})=(e_{p}\ast e_{q},\text{ }c_{a}\phi
_{pq}c_{b}) \tag{D5.1}%
\end{equation}
and $\phi_{pq}\in\Phi$ for every $e_{p},e_{q}\in E.$ If $\left(
C,\Phi\right)  $ is a mono-$\phi$ (when $\phi_{pq}=\phi$ for all $p,q$), the
block product is called a \textbf{direct product}. Two elements $(e_{i}%
,c_{x})$, $(e_{j},c_{y})\in\mathcal{L}$ are \textbf{equal}, that is,
$(e_{i},c_{x})=(e_{j},c_{y}),$ if and only if $e_{i}=e_{j}$ and $c_{x}=c_{y}.$
\end{definition}

By definition, $\Phi$ is a set of $k^{2}$ quasigroup operations $\phi_{pq}$.
The simplest case is that for which $\phi_{pq}=\phi$ for all $p,q$, that is,
if all operations of $\Phi$ are equal as in the simple direct product.
Otherwise, each system $\left(  C,\phi_{pq}\right)  $ under $\left(
C,\Phi\right)  $ can assume many forms which determine the properties of the
resulting block product $\left(  \mathcal{L},\diamond\right)  .$ The
properties of $\left(  \mathcal{L},\diamond\right)  $ will therefore depend
completely on the operation $\ast$ and the nature of the operations $\phi
_{pq}\in\Phi$. By making suitable assumptions about the systems $\left(
E,\ast\right)  $ and $\left(  C,\Phi\right)  $ we can therefore construct
composite quasigroup-type structures with desired properties. Because of this,
we shall call $\left(  E,\ast\right)  $ and $\left(  C,\Phi\right)  $ the
\textbf{\emph{generating systems}} of the block product $\left(
\mathcal{L},\diamond\right)  $.

The block product of Definition 5 satisfies Definition 3. To show this, we
partition the $n=km$ elements of $\mathcal{L}$ into $k$ cells $B_{i}%
=\{(e_{i},c_{u})\mid u=1,...,m\},$ each of order $m$. If this is done all of
the $m$ elements $(e_{i},c_{u})\in B_{i}$ will have the same $e$-component
$e_{i}$, where $i$ has a fixed value. Thus, if $(e_{p},c_{a})\in B_{p}$ and
$(e_{q},c_{b})\in B_{q}$, then by Eq. (D5.1), their product $(e_{p}\ast
e_{q},$ $c_{a}\phi_{pq}c_{b})\in B_{r},$ where $e_{p}\ast e_{q}=e_{r}$ is in
$\left(  E,\ast\right)  $. By Definition 3, we can therefore introduce the
operation $\times$ of \emph{cell multiplication} and write: $B_{p}\times
B_{q}=B_{r}.$ This simply means that given any representative element of cell
$B_{p}$ and any representative element of cell $B_{q}$, then their product is
some element of cell $B_{r}$, that is, the operation $\times$ is well defined.
The set $B=\{B_{i}\mid i=1,...,k\}$ is therefore closed under $\times$ and
hence $\left(  B,\times\right)  $ is at least a groupoid of order $k.$
Clearly, the operations $\ast$ and $\times$ are seen to be related. This is
evident from the fact that it is the operation $\ast$ of $\left(
E,\ast\right)  $ that induces the operation $\times$ on $B.$ Here we see,
however, that the operations $\phi_{pq}$ of $\Phi$ are not directly related to
$\times.$

Let $B_{p}\times B_{q}=B_{r}.$ If we form all of the $m^{2}$ binary products
of the $m$ elements of $B_{p}$ and the $m$ elements of $B_{q}$ and arrange
them in a table, they will form an $m\times m$ block of entries, denoted by
$[B_{pq}]\equiv\lbrack B_{p}\times B_{q}],$ all of which are elements of
$B_{r}$. Moreover, since there are exactly $k^{2}$ blocks $[B_{pq}]$, then
there are exactly $k^{2}$ operations $\phi_{pq}\in\Phi$ each of which is
defined as a \textbf{\emph{local operation}} over a block $[B_{pq}]$.
Henceforth, we shall use this terminology and notation for the operations of
the set $\Phi.$ Since each $\phi_{pq}\in\Phi$ is a quasigroup-type operation,
then each block $[B_{pq}]$ is a Latin square.

The idea of the block product can also be defined for groupoids in general. In
this case, the block product of two groupoids will also be a groupoid.
Moreover, there is a generalization of the direct product called a
\textbf{\emph{quasidirect product} }[6]. This can easily be shown to be a
special case of the block product.

\subsubsection{\textbf{Elementary Properties of Block Products}}

Because of the importance of block product systems in the theory of loops and
in algebra, we will now determine some of their elementary properties. In what
follows, we shall consider block products of the type given in Definition 5.

As defined, the generating systems $\left(  E,\ast\right)  $ and $\left(
C,\Phi\right)  $ of a block product $\left(  \mathcal{L},\diamond\right)  $
are required only to be at least quasigroups (axiom type \emph{A[1,4]}). This
implies that the block product $\left(  \mathcal{L},\diamond\right)  $ is also
a quasigroup as shown by:

\begin{theorem}
Let $\left(  E,\ast\right)  $ and $\left(  C,\Phi\right)  $ be quasigroups.
Then their block product $\left(  \mathcal{L},\diamond\right)  =\left(
E,\ast\right)  $\textsf{X}$\left(  C,\Phi\right)  $ is also a quasigroup, and conversely.
\end{theorem}

%

\proof
Let $(e_{p},c_{a}),$ $(e_{q},c_{b})\in\mathcal{L}.$ By Eq. (D5.1) of
Definition 5, $(e_{p},c_{a})\diamond(e_{q},c_{b})=(e_{p}\ast e_{q},c_{a}%
\phi_{pq}c_{b}).$ Since $\left(  E,\ast\right)  $ and $\left(  C,\Phi\right)
$ are quasigroups, then $e_{p}\ast e_{q}\in E,$ $c_{a}\phi_{pq}c_{b}\in C$ and
hence $(e_{p}\ast e_{q},c_{a}\phi_{pq}c_{b})\in\mathcal{L}.$ This implies that
$\left(  \mathcal{L},\diamond\right)  $ satisfies \emph{A1 }and is at least a
groupoid. Now consider the equation%

\[
(e_{p},c_{a})\diamond(e_{u},c_{x})=(e_{q},c_{b})
\]
By Eq. (D5.1), $(e_{p},c_{a})\diamond(e_{u},c_{x})=(e_{p}\ast e_{u},c_{a}%
\phi_{pu}c_{x}),$ where $e_{p}\ast e_{u}\in E$ and $c_{a}\phi_{pu}c_{x}\in C.$
Since by hypothesis $\left(  E,\ast\right)  $ and $\left(  C,\Phi\right)  $
are quasigroups, then given $e_{p},e_{q}\in E$ and $c_{a},c_{b}\in C,$ there
exist unique elements $e_{u}\in E$ and $c_{x}\in C$ such that $e_{p}\ast
e_{u}=e_{q},$ and $c_{a}\phi_{pu}c_{x}=c_{b}.$ Hence, the element
$(e_{u},c_{x})\in\mathcal{L}$ exists and is unique. Similarly, we can show
that there exists a unique element $(e_{v},c_{y})\in\mathcal{L}$ such that:%

\[
(e_{v},c_{y})\diamond(e_{p},c_{a})=(e_{q},c_{b})
\]
Therefore, $\left(  \mathcal{L},\diamond\right)  $ also satisfies \emph{A4}
from which it follows that it is at least of type \emph{A[1,4]} and is a
quasigroup. Conversely, if $\left(  \mathcal{L},\diamond\right)  $ is a
quasigroup, then $(e_{p},c_{a})\diamond(e_{q},c_{b})=(e_{p}\ast e_{q}%
,c_{a}\phi_{pq}c_{b})\in\mathcal{L}$ and is unique. This implies that
$e_{p}\ast e_{q}\in E$ and $c_{a}\phi_{pq}c_{b}\in C$ are also unique so that
$\left(  E,\ast\right)  $ and $\left(  C,\Phi\right)  $ must also be at least
of type \emph{A[1,4]} and are quasigroups. $\blacksquare$\vspace{0.05in}

The next theorem shows that there is a natural partition of $\mathcal{L}$ that
is induced by the operation $\ast$ of $\left(  E,\ast\right)  .$

\begin{theorem}
Let $\left(  \mathcal{L},\diamond\right)  =\left(  E,\ast\right)  $%
\textsf{X}$\left(  C,\Phi\right)  $ be a block product of order $n=km,$ where
$\left(  E,\ast\right)  $ and $\left(  C,\Phi\right)  $ are quasigroups,
$E=\left\{  e_{1},\ldots,e_{k}\right\}  $ and $C=\left\{  c_{1},\ldots
,c_{m}\right\}  .$ Let $B=\left\{  B_{i}\mid i=1,\ldots,k\right\}  $ where
\begin{equation}
B_{i}=\left\{  (e_{i},c_{u})\in\mathcal{L}\mid u=1,\ldots,m\right\}
\tag{T4.1}%
\end{equation}
Then (a) $B$ is a partition of $\mathcal{L}$ and (b) the set $B$ forms a
quasigroup $\left(  B,\times\right)  ,$ where $\times$ is cell multiplication.
\end{theorem}

%

\proof
(a) It follows from Eq. (T4.1) that any element $(e_{i},c_{u})\in\mathcal{L} $
is in exactly one subset $B_{i}\ $of $\mathcal{L}.$ Hence, $B_{i}\cap
B_{j}=\emptyset$ if $i\neq j$ and $\bigcup_{i=1}^{k}B_{i}=\mathcal{L}$ which
implies that $B$ is a partition of $\mathcal{L}$. (b) Now, let $(e_{i}%
,c_{u})\in B_{i}$ and let $(e_{j},c_{v})\in B_{j}.$ Then by Eq. (D5.1),
$(e_{i},c_{u})\diamond(e_{j},c_{v})=(e_{h},c_{w})\in\mathcal{L},$ where
$e_{h}=e_{i}\ast e_{j}\in E$ and $c_{w}=c_{u}\phi_{ij}c_{v}\in C.$ By Theorem
3, $\left(  \mathcal{L},\diamond\right)  $ is a quasigroup and hence
$(e_{h},c_{w})$ is unique and, by the definition of $B_{i}$ given by Eq.
(T4.1) $(e_{h},c_{w})\in B_{h},$ where $B_{h}\in B.$ This means that the
operation $\diamond$ of $\left(  \mathcal{L},\diamond\right)  $ (through the
operation $\ast$ of $\left(  E,\ast\right)  $) induces on the cells $B_{x}$ of
$B$ a well defined operation $\times$ of \emph{cell multiplication} such that
$B_{i}\times B_{j}=B_{h}$ forming a groupoid $\left(  B,\times\right)  $.
Next, consider the equation: $(e_{p},c_{u})\diamond(e_{x},c_{v})=(e_{q}%
,c_{w}),$ where $(e_{p},c_{u})\in B_{p},$ $(e_{x},c_{v})\in B_{x}$, and
$(e_{q},c_{w})\in B_{q}.$ Then we have $(e_{p},c_{u})\diamond(e_{x}%
,c_{v})=(e_{p}\ast e_{x},c_{u}\phi_{px}c_{v})=(e_{q},c_{w})$ which implies
that: $e_{p}\ast e_{x}=e_{q}$ and $c_{u}\phi_{px}c_{v}=c_{w}.$ Since $\left(
E,\ast\right)  $ and $\left(  C,\Phi\right)  $ are quasigroups, then $e_{x}\in
E$ and $c_{v}\in C$ are unique and hence the element $(e_{x},c_{v})\in B_{x}$
is also unique. This means that given any two elements $B_{p},B_{q}\in B,$
there exists a unique element $B_{x}\in B$ such that $B_{p}\times B_{x}%
=B_{q}.$ A similar argument also shows that there exists a unique element
$B_{y}\in B$ such that $B_{y}\times B_{p}=B_{q}.$ Therefore $\left(
B,\times\right)  $ is also a quasigroup. $\blacksquare$\vspace{0.05in}

Henceforth, we shall refer to the partition $B=\left\{  B_{1},\ldots
,B_{k}\right\}  $ of $\mathcal{L}$ defined in Theorem 4 as an
\textbf{\emph{E-partition}} of $\mathcal{L}.$

Although we have shown that $\left(  B,\times\right)  $ is a quasigroup, it is
not necessarily isomorphic to any subsystem of $\left(  \mathcal{L}%
,\diamond\right)  .$ Moreover, $\left(  \mathcal{L},\diamond\right)  $ does
not necessarily have any subsystems isomorphic to $\left(  E,\ast\right)
\
$and $\left(  C,\Phi\right)  .$ Nevertheless, the operation $\times$ of
$\left(  B,\times\right)  $ is clearly induced by the operation $\diamond$ of
$\left(  \mathcal{L},\diamond\right)  $ (through the operation $\ast$ of
$\left(  E,\ast\right)  $). This indicates that there is a close connection
between the structures of $\left(  \mathcal{L},\diamond\right)  $ and $\left(
B,\times\right)  .$

\subsubsection{\textbf{Coset Product Loops}}

Let us now consider an important form of the block product, called the
\emph{\textbf{coset product}} [3], that is closely related to the idea of the
group with a coset group (or quotient group) in finite group theory.

This is the interesting case when the block product $\left(  \mathcal{L}%
,\diamond\right)  =\left(  E,\ast\right)  $\textsf{X}$\left(  C,\Phi\right)  $
is a loop (like a NAFIL or a group) with a normal subloop $\left(
B_{1},\diamond\right)  ,$ where $B_{1}\subset B,$ and the system $\left(
B,\times\right)  $ is also a loop such that every $B_{i}\in B$ is a coset of
$B_{1}.$ For this particular case, we shall therefore call $B$ the
\emph{E-partition of }$\mathcal{L}$ $(\operatorname{mod}$\emph{\ }$B_{1}).$
This must be distinguished from other E-partitions of $\mathcal{L}$ in which
$(B_{1},\diamond)\in B$ is not a loop.

In constructing a coset product, it is important to be specific about the
nature of the generating systems $\left(  E,\ast\right)  $ and $\left(
C,\Phi\right)  $. By Definition 4, the multi-$\phi$ $\left(  C,\Phi\right)  $
represents several systems with a common set of elements $C.$ Therefore, each
of these systems $\left(  C,\phi_{pq}\right)  $ under $\left(  C,\Phi\right)
$ must be clearly defined. In particular, if $\left(  C,\phi_{pq}\right)  $ is
a loop, its identity element must be identified. This leads us to classify
multi-$\phi$ loops $\left(  C,\Phi\right)  $ into two major types:
\textbf{\emph{Type\ A}} when $\left(  C,\Phi\right)  $ has a common identity
element for all $\phi_{pq}\in\Phi$ and \textbf{\emph{Type\ B}} otherwise.

To fix the meaning of the coset product, we need to introduce a number of
important concepts.

\begin{definition}
Let $\left(  \mathcal{L},\diamond\right)  $ be a loop of order $n=km$ with a
subloop $\left(  B_{1},\diamond\right)  $ of order $m.$ Then $\left(
B_{1},\diamond\right)  $ is called \textbf{normal} if the set $B=\{B_{i}\mid
i=1,...,k\}$ of cosets of $B_{1}$ forms a loop $\left(  B,\times\right)  $ of
order $k$ called a \textbf{factor system}$.$
\end{definition}

This definition implies that the left and right cosets of $B_{1}$ are always
equal and that $B$ is an \emph{E-partition of }$\mathcal{L}$
$(\operatorname{mod}\,B_{1}).$ For a group (which is associative), this is a
sufficient condition for a subgroup to be normal. For a non-associative loop
(like the NAFIL), however, a subsystem whose cosets form a partition of
$\mathcal{L}$ is not necessarily normal. By Definition 6 the set $B$ of cosets
of $B_{1}$ must form a factor system $\left(  B,\times\right)  $ for $\left(
B_{1},\diamond\right)  $ to be normal. If $\left(  B_{1},\diamond\right)  $ is
normal in $\left(  \mathcal{L},\diamond\right)  $, we can formally write:
$\left(  \mathcal{L},\diamond\right)  /\left(  B_{1},\diamond\right)  =\left(
B,\times\right)  .\bigskip$

\begin{center}
$\underset{}{\underset{\text{3(A).}\left(  \left\{  \ell_{1},\ell_{2},\ell
_{3},\ell_{4}\right\}  ,\diamond\right)  \text{\ is normal}}{
\begin{tabular}
[c]{|c|cccc|cccc|}\hline
$\diamond$ & $\ell_{1}$ & $\ell_{2}$ & $\ell_{3}$ & $\ell_{4}$ & $\ell_{5}$ &
$\ell_{6}$ & $\ell_{7}$ & $\ell_{8}$\\\hline
$\ell_{1}$ & $\ell_{1}$ & $\ell_{2}$ & $\ell_{3}$ & $\ell_{4}$ & $\ell_{5}$ &
$\ell_{6}$ & $\ell_{7}$ & $\ell_{8}$\\
$\ell_{2}$ & $\ell_{2}$ & $\ell_{3}$ & $\ell_{4}$ & $\ell_{1}$ & $\ell_{6}$ &
$\ell_{7}$ & $\ell_{8}$ & $\ell_{5}$\\
$\ell_{3}$ & $\ell_{3}$ & $\ell_{4}$ & $\ell_{1}$ & $\ell_{2}$ & $\ell_{7}$ &
$\ell_{8}$ & $\ell_{5}$ & $\ell_{6}$\\
$\ell_{4}$ & $\ell_{4}$ & $\ell_{1}$ & $\ell_{2}$ & $\ell_{3}$ & $\ell_{8}$ &
$\ell_{5}$ & $\ell_{6}$ & $\ell_{7}$\\\hline
$\ell_{5}$ & $\ell_{5}$ & $\ell_{6}$ & $\ell_{7}$ & $\ell_{8}$ & $\ell_{1}$ &
$\ell_{2}$ & $\ell_{3}$ & $\ell_{4}$\\
$\ell_{6}$ & $\ell_{6}$ & $\ell_{5}$ & $\ell_{8}$ & $\ell_{7}$ & $\ell_{2}$ &
$\ell_{1}$ & $\ell_{4}$ & $\ell_{3}$\\
$\ell_{7}$ & $\ell_{7}$ & $\ell_{8}$ & $\ell_{5}$ & $\ell_{6}$ & $\ell_{3}$ &
$\ell_{4}$ & $\ell_{1}$ & $\ell_{2}$\\
$\ell_{8}$ & $\ell_{8}$ & $\ell_{7}$ & $\ell_{6}$ & $\ell_{5}$ & $\ell_{4}$ &
$\ell_{3}$ & $\ell_{2}$ & $\ell_{1}$\\\hline
\end{tabular}
}}\qquad\underset{}{\underset{\text{3(B). }\left(  \left\{  \ell_{1},\ell
_{7}\right\}  ,\diamond\right)  \;\text{is not normal}}{
\begin{tabular}
[c]{|c|cc|cc|cc|cc|}\hline
$\diamond$ & $\ell_{1}$ & $\ell_{7}$ & $\ell_{2}$ & $\ell_{8}$ & $\ell_{3}$ &
$\ell_{5}$ & $\ell_{4}$ & $\ell_{6}$\\\hline
$\ell_{1}$ & $\ell_{1}$ & $\ell_{7}$ & $\ell_{2}$ & $\ell_{8}$ & $\ell_{3}$ &
$\ell_{5}$ & $\ell_{4}$ & $\ell_{6}$\\
$\ell_{7}$ & $\ell_{7}$ & $\ell_{1}$ & $\ell_{8}$ & $\ell_{2}$ & $\ell_{5}$ &
$\ell_{3}$ & $\ell_{6}$ & $\ell_{4}$\\\hline
$\ell_{2}$ & $\ell_{2}$ & $\ell_{8}$ & $\ell_{3}$ & $\ell_{5}$ & $\ell_{4}$ &
$\ell_{6}$ & $\ell_{1}$ & $\ell_{7}$\\
$\ell_{8}$ & $\ell_{8}$ & $\ell_{2}$ & $\ell_{7}$ & $\ell_{1}$ & $\ell_{6}$ &
$\ell_{4}$ & $\ell_{5}$ & $\ell_{3}$\\\hline
$\ell_{3}$ & $\ell_{3}$ & $\ell_{5}$ & $\ell_{4}$ & $\ell_{6}$ & $\ell_{1}$ &
$\ell_{7}$ & $\ell_{2}$ & $\ell_{8}$\\
$\ell_{5}$ & $\ell_{5}$ & $\ell_{3}$ & $\ell_{6}$ & $\ell_{4}$ & $\ell_{7}$ &
$\ell_{1}$ & $\ell_{8}$ & $\ell_{2}$\\\hline
$\ell_{4}$ & $\ell_{4}$ & $\ell_{6}$ & $\ell_{1}$ & $\ell_{7}$ & $\ell_{2}$ &
$\ell_{8}$ & $\ell_{3}$ & $\ell_{5}$\\
$\ell_{6}$ & $\ell_{6}$ & $\ell_{4}$ & $\ell_{5}$ & $\ell_{3}$ & $\ell_{8}$ &
$\ell_{2}$ & $\ell_{7}$ & $\ell_{1}$\\\hline
\end{tabular}
}}$

{\small Table 3. Cayley tables of a NAFIL loop }$\left(  \mathcal{L}%
,\diamond\right)  ${\small \ of Order n = 8}$.$
\end{center}

To understand this concept of normality, consider the non-abelian NAFIL
$\left(  \mathcal{L},\diamond\right)  $ of order $n=8$ whose Cayley table is
shown in two forms in Table 3. Analysis (using the software FINITAS [4]) shows
that $\left(  \mathcal{L},\diamond\right)  $ has 8 non-trivial subsystems all
of which are groups: 3 of order 4 (1 isomorphic to $C_{4}$ and 2 are
isomorphic to $K_{4}$) and 5 of order 2 (isomorphic to $C_{2}$). Out of these
8 subgroups, 6 have cosets that form partitions of $\mathcal{L}$. In
particular, the 3 subgroups of order 4 are all \textbf{\emph{normal}}, that
is, their cosets form partitions $B$ of $\mathcal{L}$ such that $\left(
B,\times\right)  $ is a factor system of $\left(  \mathcal{L},\diamond\right)
.$ Table 3(A) shows $\left(  \mathcal{L},\diamond\right)  $ arranged in terms
of the normal subsystem $\left(  \left\{  \ell_{1},\ell_{2},\ell_{3},\ell
_{4}\right\}  ,\diamond\right)  $ of order 4 while Table 3(B) shows $\left(
\mathcal{L},\diamond\right)  $ arranged in terms of the subsystem $\left(
\left\{  \ell_{1},\ell_{7}\right\}  ,\diamond\right)  $ of order 2 which is
not normal.

The three subgroups of order 2 whose cosets also form partitions, however, are
not all normal. For instance, the subgroup $\left(  \mathcal{B}_{1}%
,\diamond\right)  \cong C_{2},$ where $\mathcal{B}_{1}=\{\ell_{1},\ell_{7}\},$
has the following cosets: $\mathcal{B}_{1}=\{\ell_{1},\ell_{7}\},$
$\mathcal{B}_{2}=\{\ell_{2},\ell_{8}\},$ $\mathcal{B}_{3}=\{\ell_{3},\ell
_{5}\},$ $\mathcal{B}_{4}=\{\ell_{4},\ell_{6}\}$ that form the partition
$\mathcal{B}=\{\mathcal{B}_{1},\mathcal{B}_{2},\mathcal{B}_{3},\mathcal{B}%
_{4}\}$ of $\mathcal{L}.$ However, this set $\mathcal{B}\ $of cosets of
$\mathcal{B}_{1}$ is not closed under cell multiplication $\times.$ To verify
this, consider the cell $\mathcal{B}_{2}=\{\ell_{2},\ell_{8}\}.$ If
$\mathcal{B}$ is closed under $\times,$ then the product of any two elements
of $\mathcal{B}_{2}$ must belong to just one cell $\mathcal{B}_{x}%
\in\mathcal{B}.$ But $\ell_{2}\diamond\ell_{8}=\ell_{5}\in\mathcal{B}_{3}$ and
$\ell_{8}\diamond\ell_{2}=\ell_{7}\in\mathcal{B}_{1}.$ Hence, $\mathcal{B}%
_{2}\times\mathcal{B}_{2}$ is not defined and therefore $\mathcal{B}$ is not
closed under $\times.$ This means that $\left(  \mathcal{B}_{1},\diamond
\right)  $ is not normal and the set $\mathcal{B}$ of cosets of $\mathcal{B}%
_{1}$ is not an E-partition of $\mathcal{L}$ (mod\emph{\ }$\mathcal{B}_{1}%
$)$.$

In the study of coset products we will have occasion to consider multi-$\phi$
systems $\left(  C,\Phi\right)  $ containing quasigroups in which the
existence of a unique identity or unique inverses are not assumed. In such
systems, we will consider weaker forms of \emph{A2} and \emph{A3} involving
two kinds of identity elements and two kinds of inverse elements [6, 8].

\begin{definition}
Let $\left(  \mathcal{L},\diamond\right)  $ be a quasigroup. Then $1^{\lambda
}\in\mathcal{L}$ is called a \textbf{\emph{left identity}} if and only if
$1^{\lambda}\star\ell=\ell$ and $1^{\rho}\in\mathcal{L}$ is called a
\textbf{\emph{right identity}} if and only if $\ell\star1^{\rho}=\ell$ for all
$\ell\in\mathcal{L}$.
\end{definition}

\begin{definition}
Let $\left(  \mathcal{L},\star\right)  $ be a loop whose identity element is
$\ell_{1}\equiv1$ and let $\ell\in\mathcal{L}$. Then $\ell^{-\lambda}$ and
$\ell^{-\rho}$ are called the \textbf{left inverse }and\textbf{\ right
inverse} of $\ell,$ respectively, if and only if $\ell^{-\lambda}\star\ell=1$
and $\ell\star\ell^{-\rho}=1$. If $\ell^{-\lambda}=\ell^{-\rho}\equiv\ell
^{-1}$ such that $\ell^{-1}\star\ell=\ell\star\ell^{-1}=1,$ then $\ell^{-1}$
is called the \textbf{two-sided inverse} or simply the \textbf{inverse} of
$\ell$.
\end{definition}

\begin{theorem}
Let $\left(  \mathcal{L},\diamond\right)  $ be a block product of order $n=km$
whose generating system $\left(  E,\ast\right)  $ is a loop whose identity is
$e_{1}$ and $\left(  C,\Phi\right)  $ is a multi-$\phi$ system such that every
$\left(  C,\phi_{ij}\right)  $ has a common identity element $c_{1}$. Then:
(\textbf{a}) $\left(  \mathcal{L},\diamond\right)  $ is also a loop whose
identity element is $(e_{1},c_{1}).$ (\textbf{b}) If $B=\left\{  B_{1}%
,\ldots,B_{k}\right\}  $ is the $E$-partition of $\mathcal{L}$, then $\left(
B_{1},\diamond\right)  $ is a subsystem of $\left(  \mathcal{L},\diamond
\right)  .$ (\textbf{c}) The system $\left(  B,\times\right)  $, where
$\times$ is cell multiplication, is a loop isomorphic to $\left(
E,\ast\right)  $ and is a factor system of $\left(  \mathcal{L},\diamond
\right)  .$ Hence, $\left(  B_{1},\diamond\right)  $ is normal.
\end{theorem}

%

\proof
(\textbf{a}) By Theorem 3, $\left(  \mathcal{L},\diamond\right)  $ is at least
a quasigroup. Let $E=\left\{  e_{1,}...,e_{m}\right\}  $ and $C=\left\{
c_{1},...,c_{k}\right\}  .$ Since $\left(  E,\ast\right)  $ is a loop, then it
has a unique identity element, say $e_{1}.$ By hypothesis $\left(
C,\Phi\right)  $ is a multi-$\phi$ system such that every $\left(  C,\phi
_{ij}\right)  $ has a common identity element $c_{1}.$ Now, let $(e_{1}%
,c_{1}),$ $(e_{i},c_{j})\in\mathcal{L}$ and form the products:%

\begin{align*}
(e_{i},c_{j})\diamond(e_{1},c_{1})  &  =(e_{i}\ast e_{1},c_{j}\phi_{i1}%
c_{1})=(e_{i},c_{j})\\
(e_{1},c_{1})\diamond(e_{i},c_{j})  &  =(e_{1}\ast e_{i},c_{1}\phi_{1j}%
c_{j})=(e_{i},c_{j})
\end{align*}
Therefore, we find that%

\[
(e_{1},c_{1})\diamond(e_{i},c_{j})=(e_{i},c_{j})\diamond(e_{1},c_{1}%
)=(e_{i},c_{j})
\]
which implies that the element $(e_{1},c_{1})\in\mathcal{L}$ is a unique
identity element under $\diamond$. Hence, the quasigroup $\left(
\mathcal{L},\diamond\right)  $ also satisfies \emph{A2 }and is therefore a
loop (type\emph{\ A[1,4,2]}).

(\textbf{b}) If $B=\{B_{1},...,B_{k}\}$ is an E-partition of $\mathcal{L},$
then $B_{1}=\{(e_{1},c_{u})\in\mathcal{L}\mid u=1,\ldots,m\}$ and thus the
identity element $(e_{1},c_{1})\in B_{1}.$ Now let $(e_{1},c_{u}),$
$(e_{1},c_{v})\in B_{1}.$ Then $(e_{1},c_{u})\diamond(e_{1},c_{v})=(e_{1}\ast
e_{1},c_{u}\phi_{11}c_{v})=(e_{1},c_{w}),$ where $c_{w}=c_{u}\phi_{11}c_{v}\in
C.$ Clearly, $(e_{1},c_{w})\in B_{1}$ and therefore $B_{1}$ is closed under
$\diamond.$ By Theorem 1, this implies that $\left(  B_{1},\diamond\right)  $
is a subsystem of $\left(  \mathcal{L},\diamond\right)  $ and is thus also a loop.

(\textbf{c}) By Theorem 4, since $B$ is an $E$-$partition$ of $\mathcal{L}$,
then $\left(  B,\times\right)  $ is at least a quasigroup of order $k$, where
$\times$ is cell multiplication. Hence, if $B_{i},B_{j}\in B,$ then we can
write: $B_{i}\times B_{j}=B_{h},$ where $B_{h}\in B$. In (a) we have shown
that $\left(  \mathcal{L},\diamond\right)  $ is a loop whose identity is
$(e_{1},c_{1})$ and that if $(e_{1},c_{u}),$\thinspace$(e_{1},c_{v})\in
B_{1},$ then their product $(e_{1},c_{u})\diamond(e_{1},c_{v})\in B_{1}\ $and
we can write: $B_{1}\times B_{1}=B_{1}$. Now, let $(e_{1},c_{u})\in B_{1} $
and $(e_{x},c_{v})\in B_{x}.$ Form the products:
\begin{align*}
(e_{1},c_{u})\diamond(e_{x},c_{v})  &  =(e_{1}\ast e_{x},\,c_{v}\phi_{1x}%
c_{u})=(e_{x},\,c_{v}\phi_{1x}c_{u})\in B_{x}\\
(e_{x},c_{v})\diamond(e_{1},c_{u})  &  =(e_{x}\ast e_{1},\,c_{v}\phi_{x1}%
c_{u})=(e_{x},\,c_{v}\phi_{x1}c_{u})\in B_{x}%
\end{align*}
This means that we can write: $B_{1}\times B_{x}=B_{x}\times B_{1}=B_{x},$
which implies that $B_{1}$ is an identity element under $\times.$ Hence, the
quasigroup $\left(  B,\times\right)  $ satisfies A2 and is therefore a loop.
To prove that $\left(  B_{1},\diamond\right)  $ is normal, we first show that
every $B_{i}\in B$ is a left/right coset of $B_{1}.$ Let $(e_{x},c_{v}%
)\in\mathcal{L}$ and $(e_{1},c_{u})\in B_{1}.$ By Theorem 4, $B_{1}=\left\{
(e_{1},c_{u})\in\mathcal{L}\mid u=1,\ldots,m\right\}  $ so that the left and
right cosets of $B_{1}$ are:
\begin{align*}
(e_{x},c_{v})\diamond B_{1}  &  =\left\{  (e_{x}\ast e_{1},\,c_{v}\phi
_{x1}c_{u})\mid x=1,\ldots,k,\text{ }u,v=1,\cdots,m\right\} \\
B_{1}\diamond(e_{x},c_{v})  &  =\left\{  (e_{1}\ast e_{x},\,c_{u}\phi
_{1x}c_{v})\mid x=1,\ldots,k,\text{ }u,v=1,\ldots,m\right\}
\end{align*}
Since $e_{1}$ is the identity of $E,$ then $e_{1}\ast e_{x}=e_{x}\ast
e_{1}=e_{x}$ and since $C$ is closed under every $\phi_{pq}\in\Phi,$ then we
find that $(e_{x},c_{v})\diamond B_{1}=B_{1}\diamond(e_{x},c_{v})=B_{x}\in B.$
Therefore, every left coset of $B_{1}$ is also a right coset and that they
form the \emph{E-partition }$B$ \emph{of }$\mathcal{L}$ $(\operatorname*{mod}%
\,B_{1})$. Clearly, the correspondence%

\[
e_{i}\leftrightarrow B_{i},\text{ }e_{j}\leftrightarrow B_{j},\text{ }%
e_{i}\ast e_{j}\leftrightarrow B_{i}\times B_{j}%
\]
defines an isomorphism between $\left(  E,\ast\right)  $ and $\left(
B,\times\right)  .$ Since, by hypothesis, $\left(  E,\ast\right)  $ is a loop,
this shows again that $\left(  B,\times\right)  $ is a loop under cell
multiplication $\times.$ Finally, since $\left(  B,\times\right)  $ is a loop,
then it is a factor system of $\left(  \mathcal{L},\diamond\right)  $.
Therefore $\left(  B_{1},\diamond\right)  $ is normal. $\blacksquare$%
\vspace{0.1in}

In Theorem 5, we showed that if the block product $\left(  \mathcal{L}%
,\diamond\right)  =\left(  E,\ast\right)  $\textsf{X}$\left(  C,\Phi\right)  $
is a loop and $B=\left\{  B_{1},\ldots,B_{k}\right\}  $ is the $E$-partition
of $\mathcal{L}$ $(\operatorname*{mod}B_{1})$, then $\left(  B_{1}%
,\diamond\right)  $ is a normal subsystem of $\left(  \mathcal{L}%
,\diamond\right)  .$ For such a system, each cell $B_{i}\in B$ is a coset of
$B_{1}$ so that we can call the operation $\times$ of cell multiplication in
$\left(  B,\times\right)  $ as \emph{coset cell multiplication}%
\textbf{\emph{\ }}and $\left(  \mathcal{L},\diamond\right)  =\left(
E,\ast\right)  $\textsf{X}$\left(  C,\Phi\right)  $ as a \textbf{\emph{coset
product}}.

\begin{definition}
A block product$\ \left(  \mathcal{L},\diamond\right)  =\left(  E,\ast\right)
$\textsf{X}$\left(  C,\Phi\right)  $ is called a \textbf{coset product}%
\emph{\ }if and only if it is a loop with a normal subsystem $\left(
B_{1},\diamond\right)  ,$ where the set $B=\left\{  B_{1},\ldots
,B_{k}\right\}  $ of cosets of $B_{1}$ is an E-partition of $\mathcal{L}$
$(\operatorname*{mod}B_{1}).$
\end{definition}

Note that Theorem 5 holds for loops in general provided that $\left(
C,\Phi\right)  $ has a common identity element (Type A). The system $\left(
B,\times\right)  $ is the loop of cosets of $B_{1}$ under the induced
operation $\times$ and is a factor system\textbf{\ }of $\left(  \mathcal{L}%
,\diamond\right)  $ modulo $B_{1}.$ If $\left(  C,\Phi\right)  $ has no common
identity element (Type B), the resulting coset product $\left(  \mathcal{L}%
,\diamond\right)  $ will still be a loop. This will be the case if $\left(
C,\Phi\right)  $ is such that every $\left(  C,\phi_{ij}\right)  $ has an
element $c_{1}$ that is at least a right identity if $j=1$ and a left identity
if $i=1.$

We also note that $\left(  B,\times\right)  $ is isomorphic to $\left(
E,\ast\right)  .$ But $\left(  \mathcal{L},\diamond\right)  $ does not
necessarily have any normal subsystem isomorphic to $\left(  E,\ast\right)  $
and hence to $\left(  B,\times\right)  .$ If $\left(  \mathcal{L}%
,\diamond\right)  $ has such a normal subsystem isomorphic to $\left(
E,\ast\right)  ,$ we shall call $\left(  \mathcal{L},\diamond\right)  $ an
\textbf{\emph{extension}} [6] of $\left(  E,\ast\right)  $ by $\left(
C,\Phi\right)  .$\vspace{0.05in}

\begin{center}%
\[
\mathbf{Coset\ Products\ with\ Type\ A\ and}\text{ }\mathbf{B}\text{ }\left(
C,\Phi\right)  \mathbf{\ Loops}%
\]

\end{center}

The next theorem deals with coset products that are \emph{invertible loops}
(like NAFIL loops and groups). This theorem considers the case when the loops
$\left(  C,\phi_{pq}\right)  $ under $\left(  C,\Phi\right)  $ have a common
identity element (Type A).

\begin{theorem}
If $\left(  E,\ast\right)  $ and $\left(  C,\Phi\right)  $ are invertible
loops, where $\left(  C,\Phi\right)  $ has a common identity element for all
$\phi_{pq}\in\Phi$ such that $\phi_{pq}=\phi_{qp},$ then their coset product
$\left(  \mathcal{L},\diamond\right)  =\left(  E,\ast\right)  $\textsf{X}%
$\left(  C,\Phi\right)  $ is also an invertible loop.
\end{theorem}

%

\proof
By Theorem 5, $\left(  \mathcal{L},\diamond\right)  $ is a loop if its
generating systems $\left(  E,\ast\right)  $ and $\left(  C,\Phi\right)  $ are
loops. Therefore, it remains to show that if $\left(  E,\ast\right)  $ and
$\left(  C,\Phi\right)  $ are invertible loops, where $\left(  C,\Phi\right)
$ has a common identity, then so is $\left(  \mathcal{L},\diamond\right)  .$
By definition, a loop is invertible if every element has a unique inverse.
Since $\left(  E,\ast\right)  $ and $\left(  C,\Phi\right)  $ are invertible,
then every element in these loops has a unique inverse. Let $(e_{p},c_{i})$
and $(e_{q},c_{j})\in\mathcal{L}$ such that $e_{p}\ast e_{q}=e_{q}\ast
e_{p}=e_{1}$ and $c_{i}\phi_{pq}c_{j}=c_{j}\phi_{qp}c_{i}=c_{1},$ where
$\phi_{pq}=\phi_{qp},$ and $e_{1}$ and $c_{1}$ are the identity elements of
$E$ and $C,$ respectively. Form the products:
\begin{align*}
(e_{p},c_{i})\diamond(e_{q},c_{j})  &  =(e_{p}\ast e_{q},\,c_{i}\phi_{pq}%
c_{j})=(e_{1},c_{1})\\
(e_{q},c_{j})\diamond(e_{p},c_{i})  &  =(e_{q}\ast e_{p},\,c_{j}\phi_{qp}%
c_{i})=(e_{1},c_{1})
\end{align*}
where $(e_{1},c_{1})$ is the identity element of $\left(  \mathcal{L}%
,\diamond\right)  .$ Thus we find that
\[
(e_{p},c_{i})\diamond(e_{q},c_{j})=(e_{q},c_{j})\diamond(e_{p},c_{i}%
)=(e_{1},c_{1})
\]
which implies that $(e_{p},c_{i})$ and $(e_{q},c_{j})$ are unique inverses in
$\left(  \mathcal{L},\diamond\right)  .$ This proves the theorem.
$\blacksquare$\vspace{0.15in}

Since groups and NAFIL loops are both invertible loops, then the coset product
in Theorem 6 is not necessarily a NAFIL; it could also be a group. However,
there are certain conditions under which the coset product is a NAFIL.

\begin{corollary}
If $\left(  E,\ast\right)  $ and $\left(  C,\Phi\right)  $ are invertible
loops, then $\left(  \mathcal{L},\diamond\right)  =\left(  E,\ast\right)
$\textsf{X}$\left(  C,\Phi\right)  $ is a NAFIL if $\left(  E,\ast\right)  $
or any of the loops $\left(  C,\phi_{pq}\right)  $ under $\left(
C,\Phi\right)  $ is a NAFIL.
\end{corollary}

This simply means that if at least one of the generating systems $\left(
E,\ast\right)  $ or $\left(  C,\Phi\right)  $ of $\left(  \mathcal{L}%
,\diamond\right)  $ is a NAFIL, then $\left(  \mathcal{L},\diamond\right)  $
will also be a NAFIL. However, $\left(  \mathcal{L},\diamond\right)  $ can
also be a NAFIL even if its generating systems are all groups. This is the
case if at least two of the groups $\left(  C,\phi_{pq}\right)  $ under
$\left(  C,\Phi\right)  $ are not isomorphic.\vspace{0.05in}

In the previous theorems, we determined what kind of coset product results if
its generating systems are given. The next theorem deals with the case when
the coset product is given as a loop and we wish to determine what its
generating systems are. Here we consider the case of a Type B loop.

\begin{theorem}
If $\left(  \mathcal{L},\diamond\right)  =\left(  E,\ast\right)  $%
\textsf{X}$\left(  C,\Phi\right)  $ is a coset product that is an invertible
loop, then:

\begin{enumerate}
\item[(\textbf{a})] $\left(  E,\ast\right)  $ is an invertible loop.

\item[(\textbf{b})] The systems under $\left(  C,\Phi\right)  $ are such that:
(\textbf{b1}) $\left(  C,\phi_{pq}\right)  $ is at least a quasigroup, where
$c_{1}$ is at least a right identity if $q=1$ and at least a left identity if
$p=1$; in particular, $\left(  C,\phi_{11}\right)  $ must be at least a loop
and (\textbf{b2}) if $p,q\neq1,$ then $\left(  C,\phi_{pq}\right)  $ and
$\left(  C,\phi_{qp}\right)  $ must be at least loops with $c_{1}$ as a common
identity element such that either $\phi_{pq}=\phi_{qp}$, in which case the
loops are invertible, or else every left (right) inverse of $\left(
C,\phi_{pq}\right)  $ is a right (left) inverse of $\left(  C,\phi
_{qp}\right)  $.
\end{enumerate}
\end{theorem}

%

\proof
(\textbf{a}) Since $\left(  \mathcal{L},\diamond\right)  $ is an invertible
loop, then it satisfies A3. Let $(e_{1},c_{1})$ be the identity and let
$(e_{x},c_{y})$ be any element of $\mathcal{L}$. Then, by A2 and Eq. (D5.1),
it follows that%

\begin{align}
(e_{1},c_{1})\diamond(e_{x},c_{y})  &  =(e_{1}\ast e_{x},\,c_{1}\phi_{1x}%
c_{y})=(e_{x},c_{y})\tag{T7.1}\\
(e_{x},c_{y})\diamond(e_{1},c_{1})  &  =(e_{x}\ast e_{1},\,c_{y}\phi_{x1}%
c_{1})=(e_{x},c_{y})\nonumber
\end{align}
which imply that
\[
e_{1}\ast e_{x}=e_{x},\text{\quad and\quad}e_{x}\ast e_{1}=e_{x}%
\]
Hence, $\left(  E,\ast\right)  $ has a unique identity element $e_{1}$ and is
therefore a loop. Now, let $(e_{p},c_{a})$ and $(e_{q},c_{b})$ be inverses in
$\left(  \mathcal{L},\diamond\right)  .$ Then, by A3 and Eq. (D5.1), it
follows that%

\begin{align}
(e_{p},c_{a})\diamond(e_{q},c_{b})  &  =(e_{p}\ast e_{q},\,c_{a}\phi_{pq}%
c_{b})=(e_{1},c_{1})\tag{T7.2}\\
(e_{q},c_{b})\diamond(e_{p},c_{a})  &  =(e_{q}\ast e_{p},\,c_{b}\phi_{qp}%
c_{a})=(e_{1},c_{1})\nonumber
\end{align}
which imply that:%

\[
e_{p}\ast e_{q}=e_{1}\text{\quad and\quad}e_{q}\ast e_{p}=e_{1}%
\]
This means that $e_{p}$ and $e_{q}$ are inverses in $\left(  E,\ast\right)  $
and therefore it follows that the loop $\left(  E,\ast\right)  $ satisfies A3
and is at least an invertible loop.

(\textbf{b}) Let $(e_{1},c_{1})$ be the identity of $\left(  \mathcal{L}%
,\diamond\right)  $. \textbf{(b1)} By Theorem\textbf{\ }3, $\left(
C,\phi_{pq}\right)  $ is at least a quasigroup. Let $(e_{x},c_{y}%
)\in\mathcal{L}$, where $x\neq1.$ By Eqs. (T7.1), we find that:
\[
c_{1}\phi_{1x}c_{y}=c_{y}\text{\quad and\quad}c_{y}\phi_{x1}c_{1}=c_{y}%
\]
which imply that $c_{1}$ is at least a \emph{left identity} in $\left(
C,\phi_{1x}\right)  $ and a \emph{right identity} in $\left(  C,\phi
_{x1}\right)  $. However, if $x=1,$ then $c_{1}$ is a unique identity and
$\left(  C,\phi_{11}\right)  $ is at least a loop. \textbf{(b2) }Let
$(e_{p},c_{a})$ and $(e_{q},c_{b})$ be inverses in $\left(  \mathcal{L}%
,\diamond\right)  .$ Since $\left(  \mathcal{L},\diamond\right)  $ satisfies
\emph{A3}, it follows from Eqs. (T7.2) that
\[
c_{a}\phi_{pq}c_{b}=c_{1}\text{\quad and\quad}c_{b}\phi_{qp}c_{a}=c_{1}%
\]
which implies that $\left(  C,\phi_{pq}\right)  $ and $\left(  C,\phi
_{qp}\right)  $ are loops with a common identity element $c_{1}$ such that:
either $\phi_{pq}=\phi_{qp}$, in which case $c_{a}$ and $c_{b}$ are inverses
so that the loops are invertible, or else $c_{a}$ is the \emph{left inverse}
of $c_{b}$ in $\left(  C,\phi_{pq}\right)  $ and $c_{a}$ is its \emph{right
inverse} in $\left(  C,\phi_{qp}\right)  $. This completes the proof of
Theorem 7. $\blacksquare$\vspace{0.05in}

In Theorem 4, we have shown that the system $\left(  B,\times\right)  $ is a
quasigroup whose operation $\times$ is induced by the operation $\ast$ of
$\left(  E,\ast\right)  $ on $B.$ In Theorem 5, we found that if the
generating systems $\left(  E,\ast\right)  $ and $\left(  C,\Phi\right)  $ of
$\left(  \mathcal{L},\diamond\right)  $ are loops and $B=\left\{  B_{1}%
,\ldots,B_{k}\right\}  $ is the E-partition of $\mathcal{L}$ $(mod$ $B_{1})$,
then $\left(  B_{1},\diamond\right)  $ is a normal subsystem of $\left(
\mathcal{L},\diamond\right)  $ and that $\left(  B,\times\right)  $ is a loop
isomorphic to $\left(  E,\ast\right)  .$ This indicates that there is a close
connection between the structures of the systems $\left(  \mathcal{L}%
,\diamond\right)  $ and $\left(  B,\times\right)  $ which we shall now
consider [6, 8].

\begin{definition}
Let $\left(  \mathcal{L},\diamond\right)  $ and $\left(  B,\times\right)  $ be
two algebraic systems. A mapping $\theta\mathbf{:}\mathcal{L\rightarrow}B$ is
a \textbf{homomorphism} from $\left(  \mathcal{L},\diamond\right)  $ to
$\left(  B,\times\right)  $ if for any $\ell_{a},\ell_{b}\in\mathcal{L}$ the
following relation holds:
\[
\theta(\ell_{a}\diamond\ell_{b})=\theta(\ell_{a})\times\theta(\ell_{b})
\]
In this case, $\left(  B,\times\right)  $ is called the \textbf{homomorphic
image} of $\left(  \mathcal{L},\diamond\right)  .$ If the map $\theta$ is
one-to-one and onto, the homomorphism is called an \textbf{isomorphism. }
\end{definition}

This definition holds for two algebraic systems with a common axiom type. In
the non-trivial case, the map\ $\theta:\mathcal{L\rightarrow}B$ is
many-to-one. For groups, the homomorphism preserves several properties of
$\left(  \mathcal{L},\diamond\right)  $ which are carried on to $\left(
B,\times\right)  $. Thus, if $\left(  \mathcal{L},\diamond\right)  $ is a
group, then $\left(  B,\times\right)  $ is also a group. For non-associative
loops in general, only the loop properties of $\left(  \mathcal{L}%
,\diamond\right)  $ are preserved. However, since the group is also an
invertible loop, then a non-associative loop (like a NAFIL) can have a factor
system that is a group. Accordingly, we call this a \textbf{\emph{loop
homomorphism}}. An isomorphism of an algebraic system (like a quasigroup or a
loop) onto itself is called an \textbf{\emph{automorphism}}. Moreover,
$\left(  B,\times\right)  $ is also called the \emph{\textbf{kernel}} of the homomorphism.

\begin{theorem}
Let $\left(  \mathcal{L},\diamond\right)  =\left(  E,\ast\right)  $%
\textsf{X}$\left(  C,\Phi\right)  $ be a coset product of order $n=km$ that is
an invertible loop, where $\left(  C,\Phi\right)  $ is of Type A. Let
$B=\{B_{p}\mid p=1,...,k\}$ (\textbf{E-partition}) and $D=\{D_{q}\mid
q=1,...,m\}$ (\textbf{C-partition}) be partitions of $\mathcal{L}$ such that
\begin{equation}
B_{p}=\{(e_{p},c_{j})\in\mathcal{L}\mid j=1,...,m\} \tag{T8.1}%
\end{equation}

\begin{equation}
D_{q}=\{(e_{i},c_{q})\in\mathcal{L}\mid i=1,\ldots,k\} \tag{T8.2}%
\end{equation}
If $(e_{1},c_{1})$ is the identity of $\left(  \mathcal{L},\diamond\right)  ,$
where $e_{1}$ is the identity of $\left(  E,\ast\right)  $ and $c_{1}$ is the
common identity element of $\left(  C,\Phi\right)  $ for all $\phi_{ij}\in
\Phi$, then: (\textbf{a}) $\left(  B_{1},\diamond\right)  $ is a subsystem of
$\left(  \mathcal{L},\diamond\right)  $ and is isomorphic to $\left(
C,\phi_{11}\right)  .$ (\textbf{b}) The system $\left(  B,\times\right)  $,
where $\times$ is coset cell multiplication, is an invertible loop isomorphic
to $\left(  E,\ast\right)  $. (\textbf{c}) If $D_{1}=\{(e_{i},c_{1}%
)\in\mathcal{L}\mid i=1,...,k\}$, then $\left(  D_{1},\diamond\right)  $ is a
subsystem of $\left(  \mathcal{L},\diamond\right)  $. Moreover, $\left(
D_{1},\diamond\right)  $ is isomorphic to $\left(  E,\ast\right)  $ and hence
also to $\left(  B,\times\right)  $. (\textbf{d}) $\left(  \mathcal{L}%
,\diamond\right)  $ and $\left(  B,\times\right)  $ are homomorphic and
$\left(  B,\times\right)  $ is a factor system of $\left(  \mathcal{L}%
,\diamond\right)  $. Hence, , $\left(  B_{1},\diamond\right)  $ is normal.
\end{theorem}

%

\proof
(\textbf{a}) By Eq. (T8.1), we have: $B_{1}=\{(e_{1},c_{j})\in\mathcal{L}\mid
j=1,...,m\}$. Hence, if $(e_{1},c_{a})$ and $(e_{1},c_{b})$ are any two
elements of $B_{1}$, then their product is $(e_{1},c_{a})\diamond(e_{1}%
,c_{b})=(e_{1},\,c_{a}\phi_{11}c_{b})\in B_{1}$, where $\phi_{11}\in\Phi$.
Hence, $\left(  B_{1},\diamond\right)  $ is a subsystem of $\left(
\mathcal{L},\diamond\right)  $ and is therefore an invertible loop. We also
find that the correspondence:
\[
(e_{1,}c_{a})\leftrightarrow c_{a},\text{ }(e_{1,}c_{b})\leftrightarrow
c_{b},\text{ \quad}(e_{1},c_{a})\diamond(e_{1},c_{b})\leftrightarrow c_{a}%
\phi_{11}c_{b}%
\]
defines an isomorphism between $\left(  B_{1},\diamond\right)  $ and $\left(
C,\phi_{11}\right)  $.

(\textbf{b}) If $(e_{p},c_{a})\in B_{p}$ and $(e_{q},c_{b})\in B_{q}$, then
$(e_{p},c_{a})\diamond(e_{q},c_{b})=$ $(e_{p}\ast e_{q},\,c_{a}\phi_{pq}%
c_{b})=(e_{r},\,c_{a}\phi_{pq}c_{b})\in B_{r}$, where $e_{r}=e_{p}\ast
e_{q}\in E$. Therefore, we can write $B_{p}\times B_{q}=B_{r}$ where
\begin{equation}
B_{p}\times B_{q}=\{(e_{p}\ast e_{q},\,c_{a}\phi_{pq}c_{b})\mid a,b=1,\ldots
,m\} \tag{T8.3}%
\end{equation}
which implies that $\left(  B,\times\right)  $ is at least a groupoid.
Moreover, the correspondence
\[
B_{p}\leftrightarrow e_{p},\text{ }B_{q}\leftrightarrow e_{q},\text{ }%
B_{p}\times B_{q}\leftrightarrow e_{p}\ast e_{q}%
\]
defines an isomorphism between $\left(  B,\times\right)  $ and $\left(
E,\ast\right)  .$ Since $\left(  \mathcal{L},\diamond\right)  $ is an
invertible loop, then by Theorem 7 so is $\left(  E,\ast\right)  .$ Moreover,
since $\left(  B,\times\right)  $ is isomorphic to $\left(  E,\ast\right)  ,$
then $\left(  B,\times\right)  $ is also an invertible loop.

(\textbf{c}) Let $(e_{p,}c_{1})$ and $(e_{q,}c_{1})$ be elements of $D_{1}$.
Then we find that $(e_{p},c_{1})\diamond(e_{q},c_{1})=(e_{p}\ast e_{q}%
,\,c_{1}\phi_{pq}c_{1})=(e_{r},\,c_{1})\in D_{1}$, where $e_{r}=e_{p}\ast
e_{q}\in E$ and $c_{1}=c_{1}\phi_{pq}c_{1}\in C$ (since $c_{1}$ is a common
identity for all $\phi_{pq}\in\Phi$). Thus, $\left(  D_{1},\diamond\right)  $
is a subsystem of $\left(  \mathcal{L},\diamond\right)  $. Clearly, the
correspondence:
\[
(e_{p},c_{1})\leftrightarrow e_{p},\text{\quad}(e_{q},c_{1})\leftrightarrow
e_{q},\text{\quad}(e_{p},c_{1})\diamond(e_{q},c_{1})\leftrightarrow e_{p}\ast
e_{q}%
\]
is an isomorphism between $\left(  D_{1},\diamond\right)  $ and $\left(
E,\ast\right)  $ and therefore $\left(  D_{1},\diamond\right)  $ is also an
invertible loop like $\left(  E,\ast\right)  $. Since $\left(  B,\times
\right)  $ is also isomorphic to $\left(  E,\ast\right)  ,$ then $\left(
B,\times\right)  $ and $\left(  D_{1},\diamond\right)  $ are also isomorphic.

(\textbf{d}) By Definition 10, $\left(  \mathcal{L},\diamond\right)  $ and
$\left(  B,\times\right)  $ are homomorphic if the map $\theta:\mathcal{L}%
\rightarrow B$ satisfies the relation: $\theta\lbrack(e_{i},c_{u}%
)\diamond(e_{j},c_{v})]=\theta(e_{i},c_{u})\times\theta(e_{j},c_{v})$ for all
$(e_{i},c_{u}),\,(e_{j},c_{v})\in\mathcal{L}.$ This is clearly satisfied by
$\left(  \mathcal{L},\diamond\right)  $ and $\left(  B,\times\right)  .$ To
see this, let $\theta(e_{i},c_{u})\rightarrow B_{i}$ and let $\theta
(e_{j},c_{v})\rightarrow B_{j}.$ Then we find that
\begin{align*}
\theta(e_{i},c_{u})\times\theta(e_{j},c_{v})  &  \rightarrow B_{i}\times
B_{j}\\
\theta\lbrack(e_{i},c_{u})\diamond(e_{j},c_{v})]  &  =\theta(e_{i}\ast
e_{j},\,c_{u}\phi_{ij}c_{v})=\theta(e_{k},c_{w})\rightarrow B_{k}%
\end{align*}
where $e_{k}=e_{i}\ast e_{j}$ and $c_{w}=c_{u}\phi_{ij}c_{v}.$ Hence, $\left(
B,\times\right)  $ is a \emph{factor system} of $\left(  \mathcal{L}%
,\diamond\right)  .$ This implies that $\left(  B_{1},\diamond\right)  $ is
normal. $\blacksquare$\vspace{0.05in}

We have seen that $\left(  \mathcal{L},\diamond\right)  $ and $\left(
B,\times\right)  $ are homomorphic and that $\left(  B_{1},\diamond\right)  $
is a normal subsystem of $\left(  \mathcal{L},\diamond\right)  .$ Moreover,
$B_{1}$ is the identity of $\left(  B,\times\right)  .$ Thus $\theta
(e_{1},c_{u})\rightarrow B_{1}$ means that the elements $(e_{1},c_{u}%
)\in\mathcal{L}$ are mapped by $\theta$ onto the identity element $B_{1}$ of
$\left(  B,\times\right)  .$ Following the usual terminology, we shall call
$\left(  B_{1},\diamond\right)  $ the \textbf{\emph{kernel}} of the
homomorphism of $\left(  \mathcal{L},\diamond\right)  $ into $\left(
B,\times\right)  $ and call $\left(  \mathcal{L},\diamond\right)  $ a
\textbf{\emph{factorable}} system. We also note that $\left(  E,\ast\right)  $
is isomorphic to the subsystem $\left(  D_{1},\diamond\right)  .$ However,
$\left(  D_{1},\diamond\right)  $ is not necessarily normal and hence it does
not necessarily form a factor system of $\left(  \mathcal{L},\diamond\right)
.$ This is because the system $(C,\Phi)$ is a multi-$\phi.$ If it happens that
$\left(  D_{1},\diamond\right)  $ is normal, then $\left(  \mathcal{L}%
,\diamond\right)  $ is called an \textbf{\emph{extension}} of $\left(
E,\ast\right)  .$ This will be the case when $\left(  \mathcal{L}%
,\diamond\right)  $ is a direct product where $(C,\Phi)$ is a mono-$\phi$.

The cells $D_{q}$ of the \emph{C-partition} $D$ defined in Theorem 8 do not
form a factor system if at least two operations of $\Phi$ are not equal. By
Definition 4, two binary operations $\phi_{pq}$ and $\phi_{rs}$ over the same
set $C$ are said to be \emph{equal}, that is $\phi_{pq}=\phi_{rs},$ if and
only if $c_{i}\phi_{pq}c_{j}=c_{i}\phi_{rs}c_{j}$ for all $c_{i},c_{j}\in C.$
This means that if $\phi_{pq}=\phi_{rs},$ then the systems $\left(
C,\phi_{pq}\right)  $ and $\left(  C,\phi_{rs}\right)  $ have identical Cayley
tables and hence are isomorphic. Because of this, we shall call them
\emph{identically isomorphic }\textbf{\emph{\ }}This distinction is important
when dealing with quasigroup-type systems whose Cayley tables are Latin
squares. Two identically isomorphic quasigroup-type systems have identical
Cayley tables. If they are simply isomorphic, their Cayley tables are not
necessarily identical but are equivalent (isotopic) Latin squares.

Theorem 8 is true even if $\left(  \mathcal{L},\diamond\right)  $ is only a
loop that is not invertible. However, if the systems under $\left(
C,\Phi\right)  $ do not have a common identity element (Type B), then $\left(
D_{1},\diamond\right)  $ will not be a subsystem of $\left(  \mathcal{L}%
,\diamond\right)  .$ In this case, $\left(  \mathcal{L},\diamond\right)  $
will not have any subsystem isomorphic to $\left(  B,\times\right)  $ and
hence to $\left(  E,\ast\right)  .$ Thus, $\left(  \mathcal{L},\diamond
\right)  $ is not an extension of $\left(  E,\ast\right)  .$

\begin{theorem}
Every loop $\left(  L,\star\right)  $ of order $n=km$ with a non-trivial
normal subloop $\left(  \mathcal{B}_{1},\star\right)  $ of order $m$ is a
coset product of the form $\left(  \mathcal{L},\diamond\right)  =\left(
E,\ast\right)  $\textsf{X}$\left(  C,\Phi\right)  ,$ where $\left(
E,\ast\right)  $ is a loop of order $k$ and $\left(  C,\Phi\right)  $ is a
multi-$\phi$ system of order $m$.
\end{theorem}

%

\proof
To prove this theorem, it is sufficient to show that the loop $\left(
L,\star\right)  $ with the normal subsystem $\left(  \mathcal{B}_{1}%
,\star\right)  \ $is isomorphic to a loop of the form $\left(  \mathcal{L}%
,\diamond\right)  =\left(  E,\ast\right)  $\textsf{X}$\left(  C,\Phi\right)
,$ where $\left(  E,\ast\right)  $ is a loop of order $k,$ $\left(
C,\Phi\right)  $ is a multi-$\phi$ system of order $m$, and the operation
$\diamond$ is defined by the composition rule (D5.1):
\[
(e_{p},c_{a})\diamond(e_{q},c_{b})=(e_{p}\mathbf{\ast}e_{q},\,c_{a}\phi
_{pq}c_{b}).
\]

By Theorem 5 a coset product loop $\left(  \mathcal{L},\diamond\right)
=\left(  E,\ast\right)  $\textsf{X}$\left(  C,\Phi\right)  $ of order $n=km$
with a normal subloop $\left(  B_{1},\diamond\right)  $ of order $k$ can be
constructed. We now claim that by a proper choice of generating systems
$\left(  E,\ast\right)  $ and $\left(  C,\Phi\right)  ,$ it is always possible
to construct such a loop $\left(  \mathcal{L},\diamond\right)  $ that is
isomorphic to the given loop $\left(  L,\star\right)  .$

By hypothesis, $\left(  L,\star\right)  $ has a normal subloop $\left(
\mathcal{B}_{1},\star\right)  $ of order $m.$ Hence, by Definition 6, the set
$\mathcal{B=}\left\{  \mathcal{B}_{1},...,\mathcal{B}_{k}\right\}  $ of cosets
of $\mathcal{B}_{1}$ forms a factor loop $(\mathcal{B},\mathtt{X})$ of order
$k,$ where $\mathtt{X}$ is coset multiplication. For the loops $\left(
L,\star\right)  $ and $\left(  \mathcal{L},\diamond\right)  $ to be
isomorphic, they must satisfy the following necessary conditions: The coset
product loop $\left(  \mathcal{L},\diamond\right)  $ must be such that
\textbf{(a)} its normal subloop $\left(  B_{1},\diamond\right)  $ is
isomorphic to $\left(  \mathcal{B}_{1},\star\right)  $ and \textbf{(b)} its
factor loop $(B,\times)$ is isomorphic to $\left(  \mathcal{B},\mathtt{X}%
\right)  .$

Let $L=\left\{  \ell_{1},...,\ell_{km}\right\}  $ and let $\mathcal{L}%
=\{(e_{x},c_{y})\mid e_{x}\in E,$ $c_{y}\in C\}$, where $E=\{e_{1}%
,...,e_{k}\}$ and $C=\{c_{1},...,c_{m}\}$. Next, let the elements of
$\mathcal{L}$ be paired one-to-one with the elements of $L$ by the rule:
\begin{equation}
(e_{x},c_{y})\leftrightarrow\ell_{h} \tag{T9.1}%
\end{equation}
where $h=m(x-1)+y$ for all $x=1,...,k$ and $y=1,...,m.$ Partition the $n$
elements of $\mathcal{L}$ into $k$ cells $B_{x}=\{(e_{x},c_{y})\mid
y=1,...,m\},$ where $x=1,...,k,$ and let $B=\{B_{1},...,B_{k}\}$ represent
this E-partition.

To satisfy condition (a), we now let $\left(  B_{1},\diamond\right)
\cong\left(  \mathcal{B}_{1},\star\right)  $ by means of the correspondence
rule (T9.1) for $x=1$ which becomes: $(e_{1},c_{y})\leftrightarrow\ell_{y}.$
Now, if $(e_{1},c_{p}),(e_{1},c_{q})\in B_{1},$ then $(e_{1},c_{p}%
)\diamond(e_{1},c_{q})=(e_{1}\ast e_{1},$ $c_{p}\phi_{11}c_{q})=(e_{1}%
,c_{r})\in B_{1},$ where $c_{r}=c_{p}\phi_{11}c_{q}.$ By the correspondence
rule, we have: $(e_{1},c_{p})\leftrightarrow\ell_{p}$ and $(e_{1}%
,c_{q})\leftrightarrow\ell_{q}.$ Hence, we must also have $(e_{1}%
,c_{p})\diamond(e_{1},c_{q})\leftrightarrow\ell_{p}\star\ell_{q}.$ If
$\ell_{p}\star\ell_{q}=\ell_{r}$ in $\left(  \mathcal{B}_{1},\star\right)  ,$
then we must set: $(e_{1},c_{r})\leftrightarrow\ell_{r}.$ It can be shown that
it is always possible to choose $\phi_{11}$ so that this correspondence is
satisfied and is an isomorphism.

To satisfy condition (b) that $(B,\times)$ is isomorphic to $\left(
\mathcal{B},\mathtt{X}\right)  $, we simply set the correspondence:
$B_{i}\leftrightarrow\mathcal{B}_{i}$ for all $i=1,...,k.$ This is possible
because $(B,\times)$ and $\left(  \mathcal{B},\mathtt{X}\right)  $ are loops
of the same order $k,$ where $B_{1}$ is the identity of $B$ and $\mathcal{B}%
_{1}$ is the identity of $\mathcal{B}$. As in case (a), if $(e_{i},c_{p})\in
B_{i}$ and $(e_{j},c_{q})\in B_{j},$ we can set the correspondence:
$(e_{i},c_{p})\leftrightarrow\ell_{p^{\prime}}$ and $(e_{i},c_{q}%
)\leftrightarrow\ell_{q^{\prime}},$ where $\ell_{p^{\prime}}\in\mathcal{B}%
_{i}$ ($p^{\prime}=m(i-1)+p$) and $\ell_{q^{\prime}}\in\mathcal{B}_{j}$
($q^{\prime}=m(j-1)+q$). Since $(e_{i},c_{p})\diamond(e_{j},c_{q})=(e_{i}\ast
e_{j},$ $c_{p}\phi_{ij}c_{j}),$ then, by this correspondence, we have:
\begin{equation}
(e_{i},c_{p})\leftrightarrow\ell_{p^{\prime}},\text{ }(e_{i},c_{q}%
)\leftrightarrow\ell_{q^{\prime}},\text{ and }(e_{i},c_{p})\diamond
(e_{j},c_{q})\leftrightarrow\ell_{p^{\prime}}\star\ell_{q^{\prime}} \tag{T9.2}%
\end{equation}
Again, it is always possible to choose $\phi_{ij}$ so that this correspondence
is satisfied. This shows that the isomorphic correspondence $B_{i}%
\leftrightarrow\mathcal{B}_{i}$ holds true in general for the operations
$\diamond$ and $\star$ and that (T9.2) defines an isomorphism.

The above arguments (a) and (b) show that a loop $\left(  L,\star\right)  $
with a normal subloop $\left(  \mathcal{B}_{1},\star\right)  $ is isomorphic
to a coset product loop of the form $\left(  \mathcal{L},\diamond\right)
=\left(  E,\ast\right)  $\textsf{X}$\left(  C,\Phi\right)  $. This proves the
theorem. $\blacksquare$\vspace{0.1in}

This theorem also holds for any loop with a non-trivial subsystem of even
order $n=2m$. In the case of NAFIL loops, the systems $\left(  B_{1}%
,\diamond\right)  ,$ $\left(  B,\times\right)  ,$ and $\left(  E,\ast\right)
$ must also be invertible. It is important to note that the NAFIL and the
group are both invertible loops (axiom type \emph{A[1,4,2,3]}). This means
that a NAFIL can have a group as a subsystem. It is therefore possible for a
factor system $\left(  B,\times\right)  $ to be a group even if $\left(
\mathcal{L},\diamond\right)  $ is a NAFIL. This is the case when $\left(
E,\ast\right)  ,$ which is isomorphic to $\left(  B,\times\right)  ,$ is a group.

Since a loop $\left(  \mathcal{L},\diamond\right)  $ may have more than one
non-trivial normal subsystem, then such a loop can be expressed as a coset
product with respect to each of its normal subsystems. This is related to the
various coset decompositions of a loop with respect to its normal subsystems:
each decomposition determines a unique E-partition of $\mathcal{L}$.\bigskip\ 

\paragraph{\textbf{Nucleus and Center of a Loop.}}

In a loop like an invertible loop, there exist certain elements (like the
identity) that associate with all elements of the loop. Such elements form
special subsets called the \textbf{\emph{nuclei}} and \emph{\textbf{center}
}of the loop [4, 5].

\begin{definition}
The \textbf{left nucleus} of a loop $\left(  \mathcal{L},\diamond\right)  $ is
the set of all elements $a\in\mathcal{L}$ such that $a(xy)=(ax)y$ for all
$x,y\in\mathcal{L}.$ The \textbf{middle nucleus} and \textbf{right nucleus}
are similarly defined. The intersection of these three nuclei is simply called
the \textbf{nucleus} $\mathcal{N}(\mathcal{L})$ of the loop. The set of all
elements in the nucleus which commute with all elements of the loop is called
the \textbf{center} $Z(\mathcal{L)}$ of the loop.
\end{definition}

Every nucleus is a subloop of the loop. It can be shown that the center is a
normal subloop of the loop and thus forms the \textbf{\emph{kernel}} of a loop
homomorphism. Clearly, every loop has a nucleus and hence also a center. For
instance, the identity element of a loop always forms a trivial center.

For a group $G$ (which is associative), the center is simply $Z(G)=\{b\in
G\mid ab=ba$ for all $a\in G\}.$ Hence, the center $Z(G)$ is always an abelian
subgroup. Moreover, it can be shown that $Z(G)$ is normal. For abelian groups,
every normal subgroup thus qualifies as a center. Since an abelian group may
have several normal subgroups, we shall only consider the intersection of all
normal subgroups of an abelian group to be its \textbf{\emph{proper center}}.%

\[
\underset{\text{Table 4(A). NAFIL loop }\left(  \mathcal{L},\diamond\right)
}{\underset{}{
\begin{tabular}
[c]{|c|cc|cc|cc|}\hline
$\diamond$ & 1 & 2 & 3 & 4 & 5 & 6\\\hline
1 & 1 & 2 & 3 & 4 & 5 & 6\\
2 & 2 & 1 & 4 & 3 & 6 & 5\\\hline
3 & 3 & 4 & 5 & 6 & 1 & 2\\
4 & 4 & 3 & 6 & 5 & 2 & 1\\\hline
5 & 5 & 6 & 1 & 2 & 4 & 3\\
6 & 6 & 5 & 2 & 1 & 3 & 4\\\hline
\end{tabular}
}}\quad\Longrightarrow\quad\underset{\text{Table 4(B). Factor group }\left(
B,\times\right)  \,\cong\,C_{3}}{\underset{}{
\begin{tabular}
[c]{|c|c|c|c|}\hline
$\times$ & $B_{1}$ & $B_{2}$ & $B_{3}$\\\hline
$B_{1}$ & $B_{1}$ & $B_{2}$ & $B_{3}$\\\hline
$B_{2}$ & $B_{2}$ & $B_{3}$ & $B_{1}$\\\hline
$B_{3}$ & $B_{3}$ & $B_{1}$ & $B_{2}$\\\hline
\end{tabular}
}}%
\]
\bigskip

For NAFIL loops in general, $Z(\mathcal{L})=\{b\in\mathcal{N}(\mathcal{L})\mid
ab=ba$ for all $a\in\mathcal{L}\},$ where $\mathcal{N}(\mathcal{L})$ is the
nucleus of $\left(  \mathcal{L},\diamond\right)  .$ This means that every
element of $Z(\mathcal{L)}$ must associate with all elements of $\mathcal{L}$.
Moreover, $Z(\mathcal{L})$ is always a normal subloop of $\left(
\mathcal{L},\diamond\right)  $. For instance, the center of the abelian NAFIL
$\left(  \mathcal{L},\diamond\right)  $ shown above is $Z(\mathcal{L}%
)=\{1,2\}$ which is a normal subgroup of $\left(  \mathcal{L},\diamond\right)
.$ It is easy to verify from the Cayley table that the elements of
$Z(\mathcal{L})$ commute and associate with all elements of $\mathcal{L}$.
Moreover, we also see that the cosets $B_{1}=\{1,2\}$, $B_{2}=\{3,4\}$,
$B_{3}=\{5,6\}$ of this center $Z(\mathcal{L})\equiv B_{1}$ form a factor
group $\left(  B,\times\right)  $ of order $m=3,$ where $B=\{B_{1},B_{2}%
,B_{3}\},$ isomorphic to the cyclic group $C_{3}.$ This also shows that
$Z(\mathcal{L})$ is normal.

\begin{theorem}
Every loop with a non-trivial center is a coset product.
\end{theorem}

%

\proof
This follows easily from Theorem 9 since the center of a loop is always a
normal subloop of the loop. $\blacksquare$\bigskip

This theorem holds for all loops like NAFIL loops and groups. Its importance
lies in the fact that the center of a loop is a very special normal subloop.
The unique decomposition of a loop in terms of the cosets of its center
enables us to determine many aspects of its structure. For instance, the
center $Z(\mathcal{L})$ of a loop $\left(  \mathcal{L},\diamond\right)  $
always forms the kernel of a loop homomorphism [6, 8] by means of which it is
possible to obtain an \textbf{\emph{upper }}or\textbf{\emph{\ ascending
central series }}[6, 8]: $1\subseteq Z_{1}\subseteq Z_{2}\subseteq
Z_{3}\subseteq...,$ where $Z_{1}=Z(\mathcal{L})$ and $Z_{i+1}$ is the full
pre-image in $\mathcal{L}$ of $Z(\mathcal{L}/Z_{i}).$\vspace{0.15in}

Other properties of coset products include:

\begin{theorem}
The coset product $\left(  \mathcal{L},\diamond\right)  =\left(
E,\ast\right)  $\textsf{X}$\left(  C,\Phi\right)  $ is abelian if
$\mathbf{\ast}$ and every $\phi_{pq}\in\Phi$ are commutative.\medskip
\end{theorem}

The proof of this theorem is trivial.\medskip

An important example of a coset product with Type B $\left(  C,\Phi\right)  $
loops is the \textbf{\emph{Octonion loop}} $\left(  \mathbb{O}_{L}%
,\star\right)  =\left(  E,\ast\right)  $\textsf{X}$\left(  C,\Phi\right)  $ of
order 16 which is a NAFIL [3]. The center of this loop is the cyclic group
$C_{2}$ of order 2 and its factor group is an abelian group of order 8
(isomorphic to the Klein group $K_{8}$). However, the octonion does not have
any abelian subgroup of order 8.\vspace{0.15in}

\subsubsection{\textbf{Direct Products}}

The block product reduces to the \emph{direct product }if $\Phi$ is a
mono-$\phi.$ Although the direct product structure is primarily associated
with groups and related associative systems, many quasigroups also have this
structure. In particular, a direct product is always an extension of $\left(
E,\ast\right)  $ by $\left(  C,\Phi\right)  .$ The following theorems will
deal with some of the basic properties of direct products in general.

\begin{theorem}
The direct product $\left(  \mathcal{L},\diamond\right)  =\left(
E,\ast\right)  $\textsf{X}$\left(  C,\Phi\right)  $ is a group if both
$\left(  E,\ast\right)  $ and $\left(  C,\Phi\right)  $ are groups, and
conversely. Otherwise, it is a NAFIL if $\left(  C,\Phi\right)  $ or $\left(
E,\ast\right)  ,$ or both are NAFIL loops.
\end{theorem}

\begin{theorem}
Let $\left(  \mathcal{L},\diamond\right)  =\left(  E,\ast\right)  $%
\textsf{X}$\left(  C,\Phi\right)  $ be a direct product of type A[1,4,2,3],
where $\left(  E,\ast\right)  $ and $\left(  C,\Phi\right)  $ have proper
subsystems $\left(  \overline{E},\mathbf{\ast}\right)  $ and $\left(
\overline{C},\Phi\right)  $, respectively. Then $\left(  \mathcal{L}%
,\diamond\right)  $ has subsystems isomorphic to the following: $\left(
E,\ast\right)  $, $\left(  \overline{E},\mathbf{\ast}\right)  $, $\left(
C,\Phi\right)  $, $\left(  \overline{C},\Phi\right)  $, $\left(
E,\ast\right)  $\textsf{X}$\left(  \overline{C},\Phi\right)  $, $\left(
\overline{E},\ast\right)  $\textsf{X}$\left(  C,\Phi\right)  ,$ and $\left(
\overline{E},\mathbf{\ast}\right)  $\textsf{X}$\left(  \overline{C}%
,\Phi\right)  $.
\end{theorem}

Any subsystem of $\left(  \mathcal{L},\diamond\right)  $ that is isomorphic to
either $\left(  E,\ast\right)  $ or $\left(  C,\Phi\right)  $ will be called
an \textbf{\emph{isomorphic image}} of $\left(  E,\ast\right)  $ or $\left(
C,\Phi\right)  $ in $\left(  \mathcal{L},\diamond\right)  $.

\begin{theorem}
Let $\left(  \mathcal{L},\diamond\right)  =\left(  E,\ast\right)  $%
\textsf{X}$\left(  C,\Phi\right)  $ be a direct product of type A[1,4,2,3]
(invertible loop) and let $\left(  \overline{E},\diamond\right)  $ and
$\left(  \overline{C},\diamond\right)  $ be the isomorphic images in $\left(
\mathcal{L},\diamond\right)  $ of $\left(  E,\ast\right)  $ and $\left(
C,\Phi\right)  $, respectively, where
\begin{align*}
\overline{E}  &  =\{(e_{i},c_{1})\in\mathcal{L}\mid i=1,\ldots,k\,\}\\
\overline{C}  &  =\{(e_{1},c_{j})\in\mathcal{L}\mid j=1,\ldots,m\}
\end{align*}
Then the set $D=\{D_{q}\mid q=1,...,m\}$ of cosets of $\overline{E}$
(\textbf{C-partition}) and the set $B=\{B_{p}\mid p=1,...,k\}$ of cosets of
$\overline{C}$ (\textbf{E-partition}) form factor systems $\left(
D,\times_{D}\right)  $ and $\left(  B,\times_{B}\right)  $ isomorphic to
$\left(  C,\Phi\right)  $ and $\left(  B,\times\right)  $, respectively, where
$\times_{B}$ and $\times_{D}$ are operations of coset cell multiplication.
\end{theorem}

This theorem shows that the direct product $\left(  \mathcal{L},\diamond
\right)  =\left(  E,\ast\right)  $\textsf{X}$\left(  C,\Phi\right)  $, which
contains subsystems isomorphic to $\left(  E,\ast\right)  $ and $\left(
C,\Phi\right)  ,$ is an \emph{extension} of $\left(  E,\ast\right)  $ by
$\left(  C,\Phi\right)  .$

\begin{theorem}
If $\left(  E,\ast\right)  $ and $\left(  C,\Phi\right)  $ are any two systems
of type A[1,4,2,3], then the direct products $\left(  E,\ast\right)
$\textsf{X}$\left(  C,\Phi\right)  $ and $\left(  C,\Phi\right)  $%
\textsf{X}$\left(  E,\ast\right)  $ are isomorphic.\vspace{0.05in}
\end{theorem}

The proofs of the above theorems are simple.

\subsection{\textbf{Non-Lagrangian Systems}}

One unique characteristic of certain NAFIL loops that distinguishes them from
groups is that they are \textbf{\emph{non-Lagrangian}}, that is, they have
non-trivial subsystems whose orders are not divisors of the order of their
underlying set. Thus, the smallest NAFIL $(L_{5},\ast)$ of order $n=5$ is
non-Lagrangian because it has subsystems of order $m=2$. In fact there are
numerous families of non-Lagrangian NAFIL loops.

For convenience, a subsystem will be called a \textbf{\emph{non-divisor}} if
its order is not a divisor of the order of its parent system. Otherwise, it
will be called a \textbf{\emph{divisor}}.

In a non-Lagrangian NAFIL, the cosets of a non-divisor subsystem do not
determine a unique partition of the underlying set as is the case with coset
groups. Hence, such subsystems are not \emph{normal}.

There are also NAFIL loops in which all non-trivial subsystems are
non-divisors. We call such NAFIL loops \textbf{\emph{anti-Lagrangian}}\emph{.}
Such a NAFIL has no non-trivial normal subsystem and is therefore simple.%

\[
\underset{\text{Table 5. Anti-Lagrangian NAFIL }\left(  L_{9},\ast\right)
\text{ of order }n\,=\,9\text{.}}{\underset{}{
\begin{tabular}
[c]{|l|llll|lllll|}\hline
$\ast$ & 1 & 2 & 3 & 4 & 5 & 6 & 7 & 8 & 9\\\hline
1 & 1 & 2 & 3 & 4 & 5 & 6 & 7 & 8 & 9\\
2 & 2 & 1 & 4 & 3 & 6 & 5 & 8 & 9 & 7\\
3 & 3 & 4 & 1 & 2 & 7 & 8 & 9 & 6 & 5\\
4 & 4 & 3 & 2 & 1 & 8 & 9 & 5 & 7 & 6\\\hline
5 & 5 & 6 & 7 & 8 & 9 & 1 & 2 & 4 & 3\\
6 & 6 & 5 & 8 & 9 & 1 & 7 & 3 & 2 & 4\\
7 & 7 & 8 & 9 & 6 & 2 & 3 & 4 & 5 & 1\\
8 & 8 & 9 & 5 & 7 & 3 & 4 & 6 & 1 & 2\\
9 & 9 & 7 & 6 & 5 & 4 & 2 & 1 & 3 & 8\\\hline
\end{tabular}
}}%
\]

As an example, the NAFIL $\left(  L_{9},\ast\right)  $ of order $n=9$ is
anti-Lagrangian because its five non-trivial subsystems are all non-divisors:
one of order 4 and four of order 2. These subsystems are: $\left(
\{1,2,3,4\},\ast\right)  $, $\left(  \{1,2\},\ast\right)  $, $\left(
\{1,3\},\ast\right)  $, $\left(  \{1,4\},\ast\right)  $, and $\left(
\{1,8\},\ast\right)  $.

The property of being non-Lagrangian only requires that the NAFIL has at least
one non-divisor subsystem. Thus, there are also non-Lagrangian NAFIL loops
with divisor subsystems. As an example, the non-Lagrangian PAP NAFIL $\left(
L_{10},\diamond\right)  $ of order $n=10$ has subsystems of order $4$
(non-divisor) and of orders $5$ and $2$ (divisors).%

\[
\underset{\text{Table 6. Non-Lagrangian NAFIL }\left(  L_{10},\diamond\right)
\text{ with divisor and non-divisor subsystems.}}{\underset{}{
\begin{tabular}
[c]{|c|rrrrr|rrrrr|}\hline
$\diamond$ & 1 & 2 & 3 & 4 & 5 & 6 & 7 & 8 & 9 & 10\\\hline
1 & 1 & 2 & 3 & 4 & 5 & 6 & 7 & 8 & 9 & 10\\
2 & 2 & 1 & 5 & 3 & 4 & 7 & 6 & 10 & 8 & 9\\
3 & 3 & 4 & 1 & 5 & 2 & 8 & 9 & 6 & 10 & 7\\
4 & 4 & 5 & 2 & 1 & 3 & 9 & 10 & 7 & 6 & 8\\
5 & 5 & 3 & 4 & 2 & 1 & 10 & 8 & 9 & 7 & 6\\\hline
6 & 6 & 7 & 8 & 9 & 10 & 1 & 2 & 3 & 4 & 5\\
7 & 7 & 6 & 10 & 8 & 9 & 2 & 1 & 5 & 3 & 4\\
8 & 8 & 9 & 6 & 10 & 7 & 3 & 4 & 1 & 5 & 2\\
9 & 9 & 10 & 7 & 6 & 8 & 4 & 5 & 2 & 1 & 3\\
10 & 10 & 8 & 9 & 7 & 6 & 5 & 3 & 4 & 2 & 1\\\hline
\end{tabular}
}}%
\]

The NAFIL $\left(  L_{10},\diamond\right)  $ is the direct product of a cyclic
group $\left(  E,\ast\right)  $ of prime order 2 and a non-Lagrangian NAFIL
$\left(  C,\Phi\right)  $ of order 5 which has 4 subsystems of order 2.
Analysis shows that $\left(  L_{10},\diamond\right)  $ has the following
subsystems: $($\{1,2,6,7\},$\diamond)$ which is of order 4 (non-divisor), as
well as $($\{1,2,3,4,5\},$\diamond)$ of order 5 and $($\{1,6\},$\diamond)$ of
order 2 (divisors). This is a special case of the following:

\begin{theorem}
The direct product of a non-composite NAFIL (or group) of prime order and a
non-Lagrangian NAFIL is a non-Lagrangian NAFIL.
\end{theorem}

%

\proof
Let $\left(  \mathcal{L},\diamond\right)  =\left(  E,\ast\right)  $%
\textsf{X}$\left(  C,\Phi\right)  $ be a direct product of two invertible
loops, where $\left(  E,\mathbf{\ast}\right)  $ is non-composite of prime
order $k$ and $\left(  C,\Phi\right)  $ is non-Lagrangian of order $m$. To
prove that $\left(  \mathcal{L},\diamond\right)  $ is non-Lagrangian, we must
show that it contains at least one subsystem $\left(  \overline{\mathcal{L}%
},\diamond\right)  $ that is a non-divisor. By Definition 5 and Theorem 12,
$\left(  \mathcal{L},\diamond\right)  $ is a NAFIL of order $n=km.$ Now let
$\left(  \overline{E},\mathbf{\ast}\right)  $ be a non-trivial subsystem of
order $k^{\prime}$ of $\left(  E,\ast\right)  $ and let$\left(  \overline
{C},\Phi\right)  $ be a non-divisor subsystem of order $m^{\prime}$ of
$\left(  C,\Phi\right)  $. By Theorem 13, $\left(  \mathcal{L},\diamond
\right)  $ has subsystems isomorphic to $\left(  \overline{E},\mathbf{\ast
}\right)  $ and $\left(  \overline{C},\Phi\right)  ,$ respectively, and also
to $\left(  \overline{E},\mathbf{\ast}\right)  $\textsf{X}$\left(
\overline{C},\Phi\right)  $ which is a subsystem of order $\mu=k^{\prime
}m^{\prime}.$ Next, let $\left(  \overline{\mathcal{L}},\diamond\right)  $ be
the subsystem isomorphic to $\left(  \overline{E},\mathbf{\ast}\right)
$\textsf{X}$\left(  \overline{C},\Phi\right)  $. Then $\left(  \overline
{\mathcal{L}},\diamond\right)  $ is a divisor of $\left(  \mathcal{L}%
,\diamond\right)  $ if and only if $\mu$ divides $n,$ that is, if $\frac
{n}{\mu}=\operatorname{integer}.$ But $n=km $ and $\mu=k^{\prime}m^{\prime}$
so that
\[
\frac{n}{\mu}=\frac{km}{k^{\prime}m^{\prime}}=\left(  \frac{k}{k^{\prime}%
}\right)  \left(  \frac{m}{m^{\prime}}\right)
\]
Because $\left(  E,\mathbf{\ast}\right)  $ is non-composite of prime order
$k$, then its only non-trivial subsystem $\left(  \overline{E},\mathbf{\ast
}\right)  \cong\left(  E,\mathbf{\ast}\right)  $ so that $k^{\prime}=k,$
$\left(  \frac{k}{k^{\prime}}\right)  =1,$ and $\frac{n}{\mu}=\frac
{m}{m^{\prime}}.$ Since, by hypothesis, $\left(  \overline{C},\Phi\right)  $
is a non-divisor of $\left(  C,\Phi\right)  $, then its order $m^{\prime}$ is
not a divisor of $m$ and hence $\frac{n}{\mu}=\frac{m}{m^{\prime}}%
\neq\operatorname{integer}.$ This means that $\mu$ does not divide $n$ which
implies that $\left(  \mathcal{L},\diamond\right)  $ contains a non-divisor
subsystem $\left(  \mathcal{L},\diamond\right)  \cong\left(  \overline
{E},\mathbf{\ast}\right)  $\textsf{X}$\left(  \overline{C},\Phi\right)  .$
Therefore, $\left(  \mathcal{L},\diamond\right)  $ is non-Lagrangian.
Following the same argument as above, we can also show that $\left(
\mathcal{L},\diamond\right)  =\left(  E,\ast\right)  $\textsf{X}$\left(
C,\Phi\right)  $ is non-Lagrangian if $\left(  C,\Phi\right)  $ is
non-composite of prime order and $\left(  E,\ast\right)  $ is non-Lagrangian.
This completes the proof of this Theorem. $\blacksquare$

Another example of this Theorem is the direct product of the composite NAFIL
$\left(  L_{5},\ast\right)  $ of prime order 5 and the cyclic group $C_{3}$ of
order 3. This direct product is a NAFIL of order $n=15$ with non-trivial
subsystems of orders 2 and 6 (non-divisors) and 3 and 5 (divisors).

\begin{corollary}
The direct product of two composite NAFIL loops of prime order is non-Lagrangian.
\end{corollary}

%

\proof
Let $\left(  \mathcal{L},\diamond\right)  =\left(  E,\ast\right)  $%
\textsf{X}$\left(  C,\Phi\right)  $ be a direct product of two composite
invertible loops, where $\left(  E,\ast\right)  $ is of prime order $k$ and
$\left(  C,\Phi\right)  $ is of prime order $m$. Hence, $\left(
\mathcal{L},\diamond\right)  $ is an invertible loop of order $n=km.$ Now let
$\left(  \overline{E},\mathbf{\ast}\right)  $ of order $k^{\prime}$ and
$\left(  \overline{C},\Phi\right)  $ of order $m^{\prime}$ be non-trivial
subsystems of $\left(  E,\ast\right)  $ and $\left(  C,\Phi\right)  ,$
respectively. By Theorem 14, $\left(  \mathcal{L},\diamond\right)  $ has a
subsystem $\left(  \overline{\mathcal{L}},\diamond\right)  \cong\left(
\overline{E},\mathbf{\ast}\right)  $\textsf{X}$\left(  \overline{C}%
,\Phi\right)  $ which is of order $\mu=k^{\prime}m^{\prime}.$ Then $\left(
\overline{\mathcal{L}},\diamond\right)  $ is a divisor of $\left(
\mathcal{L},\diamond\right)  $ if and only if $\mu$ divides $n,$ that is, if
$\frac{n}{\mu}=\operatorname{integer}.$ But $n=km $ and $\mu=k^{\prime
}m^{\prime}$ so that $\frac{n}{\mu}=\frac{km}{k^{\prime}m^{\prime}}.$ Since
$k$ and $m$ are prime, then $\mu$ divides $m $ if and only if $k^{\prime}=m$
and $m^{\prime}=k,$ that is, if $\mu=n,$ in which case $\left(  \overline
{\mathcal{L}},\diamond\right)  \cong\left(  \mathcal{L},\diamond\right)  $.
For if $k^{\prime}=m$ and $m^{\prime}\neq k,$ then $\frac{n}{\mu}=\frac
{m}{k^{\prime}}\neq\operatorname{integer}$ since $m$ is prime. Similarly, if
$m^{\prime}=k$ and $k^{\prime}\neq m,$ then $\frac{n}{\mu}=\frac{k}{m^{\prime
}}\neq\operatorname{integer}.$ In both cases, we find that $\mu$ does not
divide $n.$ This means that $\left(  \mathcal{L},\diamond\right)  $ contains
at least one non-divisor subsystem $\left(  \overline{\mathcal{L}}%
,\diamond\right)  $ and is therefore non-Lagrangian. $\blacksquare$

The simplest example of a NAFIL satisfying this Corollary is the direct
product of the composite NAFIL $\left(  L_{5},\ast\right)  $ of prime order 5
by itself. This is a NAFIL of order $n=25$ with non-trivial subsystems of
orders 2, 4, and 10 (non-divisors) and of order 5 (divisor).

In finite group theory, a group is called \textsl{simple} if it has no proper
non-trivial normal subgroups. For instance, any group of prime order is simple
because it has no non-trivial subgroup of any kind. This idea of a system
being simple can also be extended to NAFIL loops.

\begin{definition}
An invertible loop is called \textbf{simple} if it has no proper non-trivial
normal subsystem. If it has no non-trivial subsystem of any kind, it is called
\textbf{non-composite} or \textbf{plain }[10].
\end{definition}

In general, any invertible loop of prime order is simple. For instance, the
smallest NAFIL loop $(L_{5},\ast)$ of prime order 5 is simple. Most NAFIL
loops of prime order $n\geq7$ are simple and plain. Thus, of the 2,317
non-isomorphic abelian NAFIL loops of order 7, exactly 638 are plain [2].

\begin{center}%
\[
\underset{\text{Table 7(A). }(L_{7},\ast\ast\text{Composite}}{\underset{}{
\begin{tabular}
[c]{|c|ccc|cccc|}\hline
$\circ$ & 1 & 2 & 3 & 4 & 5 & 6 & 7\\\hline
1 & 1 & 2 & 3 & 4 & 5 & 6 & 7\\
2 & 2 & 3 & 1 & 5 & 6 & 7 & 4\\
3 & 3 & 1 & 2 & 7 & 4 & 5 & 6\\\hline
4 & 4 & 5 & 6 & 1 & 7 & 3 & 2\\
5 & 5 & 6 & 7 & 2 & 1 & 4 & 3\\
6 & 6 & 7 & 4 & 3 & 2 & 1 & 5\\
7 & 7 & 4 & 5 & 6 & 3 & 2 & 1\\\hline
\end{tabular}
}}\qquad\quad\underset{\text{Table 7(B). Non-composite}}{\underset{}{
\begin{tabular}
[c]{|c|ccccccc|}\hline
$\star$ & 1 & 2 & 3 & 4 & 5 & 6 & 7\\\hline
1 & 1 & 2 & 3 & 4 & 5 & 6 & 7\\
2 & 2 & 3 & 1 & 5 & 4 & 7 & 6\\
3 & 3 & 1 & 4 & 6 & 7 & 2 & 5\\
4 & 4 & 5 & 6 & 7 & 2 & 1 & 3\\
5 & 5 & 4 & 7 & 2 & 6 & 3 & 1\\
6 & 6 & 7 & 2 & 1 & 3 & 5 & 4\\
7 & 7 & 6 & 5 & 3 & 1 & 4 & 2\\\hline
\end{tabular}
}}%
\]

\end{center}

Tables 7 show the Cayley tables of two simple NAFIL loops\ of order
$n=7.$\ 7(A) is the Cayley table of a simple NAFIL that is composite; it has
non-trivial subsystems of orders $m=2,3$\ all of which are groups. On the
other hand, 7(B) shows that of a simple NAFIL that is non-composite
(\emph{plain}).

The identification and study of simple NAFIL loops is very important to the
understanding of NAFIL structure. As in finite group theory where simple
groups are known to be the building blocks of composite groups, simple NAFIL
loops also play a similar role. Therefore, the classification of simple NAFIL
loops is a central problem of NAFIL theory.

\begin{remark}
\emph{An invertible loop of even order }$n=2m$\emph{\ is }\textbf{\emph{not
simple}}\emph{\ if it contains at least one subsystem of order }%
$m.$\emph{\ This is true for loops in general and it follows easily from the
fact that a subsystem of order }$m=n/2$\emph{\ is always normal.}
\end{remark}

\section{Non-Composite NAFIL Loops}

The study of non-composite NAFIL loops has so far not been given much
attention by loop theorists. Very little is known about these loops called
\textbf{\emph{plain }}[7] which, in some ways, are analogous to groups of
prime order. Therefore, the identification and characterization of plain NAFIL
loops is an important problem of NAFIL theory.

\subsection{\textbf{Plain NAFIL Loops}}

As in the case of prime order groups, the simplest kind of simple NAFIL is one
that is not-composite: it has no non-trivial subsystems of any kind. There are
many examples of this type of simple NAFIL called \textbf{\emph{plain }}to
distinguish them from those that are simple but composite. Such a system is
also called \textbf{\emph{anti-associative}} because it contains no
non-trivial groups and hence A6 is not satisfied in any non-trivial way within it.

Plain NAFIL loops are analogous to groups of prime order and therefore they
play an important role in the study of NAFIL structure. Studies have so far
shown that all NAFIL loops of order $n=5$ and $6$ are composite and that there
exists a number of plain NAFIL loops of order $n=7$ one of which is shown in
Table 7(B). Thus, the smallest plain NAFIL is of order $n=7.$ All plain NAFIL
loops are of odd order and are anti-associative.

Like the groups of prime order which constitute a family of plain groups,
there are also known families of plain NAFIL loops. One such family [7]
consists of abelian NAFIL loops of odd order $n=2m+1,\,$where $m\geq3.$
\medskip

\begin{center}
$\underset{}{\underset{\text{PLAIN 7N}}{
\begin{tabular}
[c]{|c|ccccccc|}\hline
{\small {*}} & {\small 1} & {\small 2} & {\small 3} & {\small 4} & {\small 5}
& {\small 6} & {\small 7}\\\hline
{\small 1} & {\small 1} & {\small 2} & {\small 3} & {\small 4} & {\small 5} &
{\small 6} & {\small 7}\\
{\small 2} & {\small 2} & {\small 3} & {\small 4} & {\small 5} & {\small 6} &
{\small 7} & {\small 1}\\
{\small 3} & {\small 3} & {\small 4} & {\small 2} & {\small 6} & {\small 7} &
{\small 1} & {\small 5}\\
{\small 4} & {\small 4} & {\small 6} & {\small 5} & {\small 7} & {\small 1} &
{\small 3} & {\small 2}\\
{\small 5} & {\small 5} & {\small 6} & {\small 7} & {\small 1} & {\small 4} &
{\small 2} & {\small 3}\\
{\small 6} & {\small 6} & {\small 7} & {\small 1} & {\small 2} & {\small 3} &
{\small 5} & {\small 4}\\
{\small 7} & {\small 7} & {\small 1} & {\small 5} & {\small 3} & {\small 2} &
{\small 4} & {\small 6}\\\hline
\end{tabular}
}}\quad\underleftrightarrow{transpose}\quad\underset{}{\underset{\text{PLAIN
7N-T}}{
\begin{tabular}
[c]{|c|ccccccc|}\hline
{\small {*}'} & {\small 1} & {\small 2} & {\small 3} & {\small 4} & {\small 5}
& {\small 6} & {\small 7}\\\hline
{\small 1} & {\small 1} & {\small 2} & {\small 3} & {\small 4} & {\small 5} &
{\small 6} & {\small 7}\\
{\small 2} & {\small 2} & {\small 3} & {\small 4} & {\small 5} & {\small 6} &
{\small 7} & {\small 1}\\
{\small 3} & {\small 3} & {\small 4} & {\small 2} & {\small 6} & {\small 7} &
{\small 1} & {\small 5}\\
{\small 4} & {\small 4} & {\small 5} & {\small 6} & {\small 7} & {\small 1} &
{\small 2} & {\small 3}\\
{\small 5} & {\small 5} & {\small 6} & {\small 7} & {\small 1} & {\small 4} &
{\small 3} & {\small 2}\\
{\small 6} & {\small 6} & {\small 7} & {\small 1} & {\small 3} & {\small 2} &
{\small 5} & {\small 4}\\
{\small 7} & {\small 7} & {\small 1} & {\small 5} & {\small 2} & {\small 3} &
{\small 4} & {\small 6}\\\hline
\end{tabular}
}}$

$\underset{}{\text{Table 8. Cayley tables of two plain NAFIL loops}}$ of order
$n=7.$
\end{center}

\begin{remark}
So far, we have discussed some of the basic properties of composite and
non-composite NAFIL loops. The important question is: How do we construct and
analyze a particular loop? In our studies, we made use of a powerful computer
software called \textbf{FINITAS }[4]. This software was developed by a team of
PUP researchers and students at the PUP SciTech R\&D Center with the support
of the Department of Science \& Technology.\bigskip
\end{remark}

\vspace{0.5in}

\newpage


\begin{thebibliography}{99}                                                                                               %


\bibitem {ref1}G. Birkhoff and S. MacLane, \emph{A Survey of Modern Algebra},
Macmilland Company, New York (1958)

\bibitem {ref2}R. E. Cawagas, \emph{Determination and Characterization of
NAFIL Loops of Small Order 5, 6, 7}, PUP Journal of Research \& Exposition,
Vol. 3, No. 1 (2004).

\bibitem {ref3}------, \emph{Finite Invertible Loops of the Coset Product
Type}, PUP Journal of Research \& Exposition (2002), Vol. 2, Issue 1, pp. 11-20.

\bibitem {ref4}------, \emph{FINITAS -- A Software for the Construction and
Analysis of Finite Algebraic Structures}, PUP Journal of Research and
Exposition, Vol. 1, No. 1 , pp. 1-10 (1997). [See also: R. E. Cawagas,
\emph{AXIOMS - Software for the Construction and Analysis of Finite
Quasigroups, Semigroups and Related Structures,} Proceedings of the Second
Asian Mathematical Conference 1995, World Scientific, Singapore-New
Jersey-London-Hong Kong (1998), 401 405.]

\bibitem {ref5}------, \emph{Introduction to Non-Associative Finite Invertible
Loops}, PUP Journal of Science \& Technology, Vol. 1, No. 2 (2007)

\bibitem {ref6}O. Chein, et al (Editors) \emph{Quasigroups and Loops: Theory
and Applications}, Sigma Series in Pure Mathematics, Helderman Verlag Berlin (1990).

\bibitem {ref7}H. Griffin, \emph{The abelian Quasi-Group}, American Journal of
Mathematics, Vol. 62, Issue 1/4 (1940), pp. 725-737.

\bibitem {ref8}H. O. Pflugfelder, \emph{Quasigroups and Loops: Introduction},
Sigma Series in Pure Mathematics, Helderman Verlag Berlin (1990).

\bibitem {ref9}D. J. S. Robinson, \emph{A Course in the Theory of Groups},
Springer Verlag, New York (1982) [See also: B. Baumslag and B. Chandler,
\emph{Group Theory}, Schaum's Outline Series in Mathematics, McGraw-Hill Book
Company (19680]

\bibitem {ref10}J. D. H. Smith, \emph{Mal'cev Varieties}, Lecture Notes in
Mathematics 554, Springer-Verlag, Berlin$\cdot$Heidelberg$\cdot$New York
(1976), pp. 96-112..
\end{thebibliography}
\end{document}